\documentclass[a4paper,12pt]{article}

\usepackage{amsfonts}
\usepackage{amsthm}
\usepackage{amssymb}
\usepackage{amsmath}
\usepackage{upref}
\usepackage{mathrsfs}
\usepackage{color}
\usepackage[citecolor=blue,colorlinks=true]{hyperref}
\usepackage{enumitem}

\newtheorem{theorem}{Theorem}
\newtheorem{proposition}[theorem]{Proposition}
\newtheorem{lemma}[theorem]{Lemma}
\newtheorem{corollary}[theorem]{Corollary}
\newtheorem{definition}[theorem]{Definition}
\newtheorem{remark}[theorem]{Remark}

\def\N{\mathbb{N}}

\def\R{\mathbb{R}}

\def\ds{\displaystyle} 
\def\div{{\rm div}}
\def\refe#1{(\ref{#1})}

\def\ocirc#1{\ifmmode\setbox0=\hbox{$#1$}\dimen0=\ht0
    \advance\dimen0 by1pt\rlap{\hbox to\wd0{\hss\raise\dimen0
    \hbox{\hskip.2em$\scriptscriptstyle\circ$}\hss}}#1\else
    {\accent"17 #1}\fi} 
    
\def\qed{\rule{0.2cm}{0.2cm}}

\def\eps{\varepsilon}

\def\<{\langle}
\def\>{\rangle}
\def\F{\mathcal{F}}
\def\GG{\mathbf{G}}
\def\P{\mathbb{P}}
\def\E{\mathbb{E}}
\def\T{\mathbb{T}}

\DeclareMathOperator{\esssup}{ess\,sup}

\begin{document}

\title{Scalar conservation laws with stochastic forcing, revised version}
\author{A. Debussche\thanks{IRMAR, ENS Cachan Bretagne, CNRS, UEB. av Robert Schuman, F-35170 Bruz, France. Email: arnaud.debussche@bretagne.ens-cachan.fr} and J. Vovelle\thanks{Universit\'e de Lyon ; CNRS ; Universit\'e Lyon 1, Institut Camille Jordan,  43 boulevard du 11 novembre 1918, F-69622 Villeurbanne Cedex, France. Email: vovelle@math.univ-lyon1.fr}}
\maketitle

\begin{abstract} We show that the Cauchy Problem for a randomly forced, periodic multi-di\-men\-sio\-nal scalar first-order conservation law with additive or multiplicative noise is well-posed: it admits a unique solution, cha\-rac\-te\-ri\-zed by a kinetic formulation of the problem, which is the limit of the solution of the stochastic parabolic approximation.
\end{abstract}

\medskip


\section{Introduction}


Let $(\Omega,\F,\P,(\F_t),(\beta_k(t)))$ be a stochastic basis and let $T>0$. In this paper, we study the first-order scalar conservation law with stochastic forcing
\begin{equation}\label{stoSCL}
du+\div(A(u))dt=\Phi(u) dW(t),\quad x\in\T^N, t\in(0,T).
\end{equation}
The equation is periodic in the space variable $x$:  $x\in\T^N$ where $\T^N$ is the $N$-dimensional torus. The flux function $A$ in \refe{stoSCL} is supposed to be of class $C^2$: $A\in C^2(\R;\R^N)$ and its derivatives have at most polynomial growth. We assume that $W$ is a cylindrical Wiener process: $W=\sum_{k\geq 1}\beta_k e_k$, where the $\beta_k$ are independent brownian processes and $(e_k)_{k\geq 1}$ is a complete orthonormal system in a Hilbert space $H$. For each $u\in\R$, $\Phi(u)\colon H\to L^2(\T^N)$ is defined by $\Phi(u)e_k=g_k(u)$ where $g_k(\cdot,u)$ is a regular function on $\T^N$. More precisely, we assume
$g_k\in C(\T^N\times\R)$, with the bounds 
\begin{align}
\GG^2(x,u)=\sum_{k\geq 1}|g_k(x,u)|^2\leq D_0(1+|u|^2),\label{D0}\\
\sum_{k\geq 1}|g_k(x,u)-g_k(y,v)|^2\leq D_1(|x-y|^2+|u-v|h(|u-v|)),\label{D1}
\end{align}
where $x,y\in\T^N$, $u,v\in\R$, and $h$ is a continuous non-decreasing function on $\R_+$ with $h(0)=0$. Note in particular that, for each $u\in\R$, $\Phi(u)\colon H\to L^2(\T^N)$ is Hilbert-Schmidt since $\|g_k(\cdot,u)\|_{L^2(\T^N)}\leq\|g_k(\cdot,u)\|_{C(\T^N)}$ and thus
\begin{equation*}
\sum_{k\geq 1}\|g_k(\cdot,u)\|_{L^2(\T^N)}^2\leq D_0(1+|u|^2).
\end{equation*}

\bigskip

The Cauchy Problem, {\it resp.} the Cauchy-Dirichlet Problem, for the stochastic equation~\refe{stoSCL} in the case of an additive noise ($\Phi$ independent on $u$) has been studied in \cite{Kim03}, {\it resp.} \cite{ValletWittbold09}. Existence and uniqueness of entropy solutions are proved in both papers. 
The Cauchy Problem for the stochastic equation~\refe{stoSCL} in case where the noise is multiplicative (and satisfies \refe{D0}-\refe{D1} above) has been studied in \cite{FengNualart08}. In  \cite{FengNualart08}, uniqueness of (strong) entropy solution is proved in any dimension, existence in dimension $1$. 
\medskip

Our purpose here is to solve the Cauchy Problem for \refe{stoSCL} in any dimension. To that purpose, we use a notion of kinetic solution, as developed by  Lions, Perthame, Tadmor for deterministic first-order scalar conservation laws~\cite{LionsPerthameTadmor94}. A very basic reason to this approach is the fact that no pathwise $L^\infty$ a priori estimates are known for \refe{stoSCL}. Thus, viewing \refe{stoSCL} as an extension of the deterministic first-order conservation law, we have to turn to the $L^1$ theory developed for the latter, for which the kinetic formulation, once conveniently adapted, is slightly better suited than the renormalized-entropy formulation (developed in \cite{CarrilloWittbold99} for example). 
\medskip

There is also a definite technical advantage to the kinetic approach, for it allows to keep track of the dissipation of the noise by solutions. For entropy solutions, part of this information is lost and has to be recovered at some stage (otherwise, the classical approach \`a la Kruzhkov~\cite{Kruzhkov70} to Comparison Theorem fails): accordingly, Feng and Nualart need to introduce a notion of ``strong" entropy solution, {\it i.e.} entropy solution satisfying the extra property that is precisely lacking~\cite{FengNualart08}. This technical difference between the notions of kinetic and entropy solution already appears in the context of degenerate parabolic equations: in the comparison of entropy solutions for such hyperbolic-parabolic equations, it is necessary to recover in a preliminary step the quantitative entropy dissipation due to the second-order part in non-degeneracy zones (see Lemma~1 in \cite{Carrillo99}). For kinetic solutions, this preliminary step is unnecessary 
 since this dissipation is already encoded in the structure of the kinetic measure, (see Definition~2.2 in \cite{ChenPerthame03}).

\medskip

In the case of an additive noise, Kim~\cite{Kim03} and Vallet and Wittbold~\cite{ValletWittbold09} introduce the auxiliary unknown $w:=u-\Phi W$ that satisfies the first-order scalar conservation law
\begin{equation}
\partial_t w+\div(B(x,t,w))=0,
\label{addSCLw}\end{equation}
where the flux $B(x,t,w):=A(w+\Phi(x)W(t))$ is non-autonomous and has limited (pathwise H\"older-) regularity with respect to the variable $t$. Then entropy solutions are defined on the basis of \refe{addSCLw}. In this way it is actually possible to avoid the use of It\^o stochastic calculus. 
\smallskip

In the case of an equation with a multiplicative noise, Feng and Nualart define a notion of entropy solution by use of regular entropies and It\^o Formula~\cite{FengNualart08}. They also define a notion of {\it strong} entropy solution, which is an entropy solution sa\-tis\-fy\-ing an additional technical criterion. This additional criterion is required to prove a comparison result between entropy and strong entropy solution. As already mentioned, they are able to prove existence of strong entropy solutions
only in dimension one. 
\smallskip

In all three papers~\cite{Kim03,FengNualart08,ValletWittbold09}, existence is proved {\it via} approximation by stochastic parabolic equation. We will proceed similarly, {\it cf.} Theorem~\ref{th:cvexists}. Consequently, our notion of solution, defined in Definition~\ref{defkineticsol}, happen to be equivalent to the notion of entropy solution used in \cite{Kim03,FengNualart08,ValletWittbold09}, provided the convergence of the vanishing viscosity method has been proved, hence in the context of \cite{Kim03,ValletWittbold09} or in \cite{FengNualart08} in dimension $1$\footnote{note that we consider periodic boundary conditions here, unlike \cite{Kim03,FengNualart08,ValletWittbold09}. However, our results extend to the whole Cauchy Problem or to the Cauchy-Dirichlet Problem.}. In fact, we prove that our notion of kinetic solution is also equivalent to the notion of (mere -- not strong) entropy solution of \cite{FengNualart08}, whatever the dimension, see section \ref{s3.3}. 
\bigskip 

Our main results states that under assumptions \refe{D0}  and \refe{D1}, there exists a unique kinetic solution in any space dimension. Due to the equivalence with entropy solution, we fill the gap 
left open in \cite{FengNualart08}. Moreover, the use of kinetic formulation considerably 
simplifies the arguments. For instance, to construct a solution, only weak compactness of 
the viscous solutions is necessary. 

\bigskip

There are related problems to \refe{stoSCL}. We refer to the references given in, {\it e.g.} \cite{Kim03,ValletWittbold09}, in particular concerning the study of the deterministic inviscid Burgers equation with random initial datum. One of the important question in the analysis of \refe{stoSCL} (and, more precisely, in the analysis of the evolution of the {\it law} of the solution process $u(t)$) is also the existence (uniqueness, ergodic character, etc.) of an invariant measure. This question has been fully addressed in \cite{EKMS00} for the inviscid periodic Burgers equation in dimension $1$ by use of the Hopf-Lax formula. 
\medskip

Our analysis of \refe{stoSCL} uses the tools developed over the past thirty years for the analysis of deterministic first-order scalar conservation laws, in particular the notion of generalized solution. Thus, in Section~\ref{sec:kisol}, we introduce the notion of solution to \refe{stoSCL} by use of the kinetic formulation, and complement it with a notion of {\it generalized} solution. In Section~\ref{sec:comparison}, we prove Theorem~\ref{th:Uadd}, which gives uniqueness (and comparison results) for solutions and also shows that a generalized solution is actually necessarily a solution. This result is used in Section~\ref{sec:exists}: we study the parabolic approximation to \refe{stoSCL} and show that it converges to a generalized solution, hence to a solution. This gives existence and uniqueness of a solution, Theorem~\ref{th:cvexists}.

\bigskip

{\bf Note:} This is an improved version of the article entitled " Scalar conservation laws with stochastic forcing" published in the Journal of Functional Analysis, 259 (2010), pp.
1014-1042. 

Since the publication of this paper, other articles on this subject have appeared. Chen, Ding, Karlsen \cite{ChenDingKarlsen} and Bauzet, Vallet and 
Wittbolt \cite{BVW12} have generalized the Kruzkov approach to the stochastic case for an equation
similar to the one treated here. Hofmanov\'a has 
proved convergence of a BGK approximation (\cite{H1}). Debussche, Hofmanov\'a and Vovelle
have treated the degenerate parabolic quasilinear case 
(\cite{DHV13}).

Also Lions, Perthame, Souganidis (\cite{LPS13}) have treated the case of a stochastic conservation law with the
stochastic term in the flux. The methods are completely different in this paper. 

\bigskip

{\bf Acknowledgement:} The authors warmly thank Martina Hofmanov\'a and Sylvain Dotti for a careful reading of our manuscript. They raised several imprecision and mistakes. These have been corrected in the present version of our work.

\section{Kinetic solution}\label{sec:kisol}


\subsection{Definition}

\begin{definition}[Kinetic measure] We say that a map $m$ from $\Omega$ to the set of non-negative finite measures over $\T^N\times[0,T]\times\R$ is a kinetic measure if
\begin{enumerate}
\item $m$ is measurable, in the sense that for each $\phi\in C_b(\T^N\times[0,T]\times\R)$, $\<m,\phi\>\colon\Omega\to\R$ is,
\item $m$ vanishes for large $\xi$: if $B_R^c=\{\xi\in\R,|\xi|\geq R\}$, then
\begin{equation}
\lim_{R\to+\infty}\E m(\T^N\times[0,T]\times B_R^c)=0,
\label{inftym}\end{equation}
\item for all $\phi\in C_b(\T^N\times\R)$, the process
\begin{equation*}
t\mapsto \int_{\T^N\times[0,t]\times\R} \phi(x,\xi)dm(x,s,\xi)
\end{equation*}
is predictable.
\end{enumerate}
\label{def:kineticmeasure}\end{definition}

\begin{definition}[Solution] Let $u_{0}\in L^\infty(\T^N)$. A measurable function $u\colon\T^N\times [0,T]\times\Omega\to\R$ is said to be a solution to~\refe{stoSCL} with initial datum $u_0$ if $(u(t))$ is predictable, for all $p\geq 1$, there exists $C_p\geq 0$ such that
\begin{equation}
\E\left(\esssup_{t\in[0,T]}\|u(t)\|_{L^p(\T^N)}^p\right)\leq C_p,
\label{eq:integrabilityu}\end{equation} 
and if there exists a kinetic measure $m$ such that $f:=\mathbf{1}_{u>\xi}$ satisfies:  for all $\varphi\in C^1_c(\T^N\times[0,T)\times\R)$, 
\begin{multline}
\int_0^T\<f(t),\partial_t \varphi(t)\>dt+\<f_0,\varphi(0)\>
+\int_0^T \<f(t),a(\xi)\cdot\nabla\varphi(t)\>dt\\
=-\sum_{k\geq 1}\int_0^T\int_{\T^N}g_k(x,u(x,t))\varphi(x,t,u(x,t)) dxd\beta_k(t)\\
-\frac{1}{2}\int_0^T\int_{\T^N} \partial_\xi\varphi(x,t,u(x,t))\GG^2(x,u(x,t)) dx dt+m(\partial_\xi\varphi),
\label{eq:kineticupre}\end{multline}
a.s., where $\GG^2:=\sum_{k=1}^\infty |g_k|^2$ and $a(\xi):=A'(\xi)$.
\label{defkineticsol}\end{definition}

In \refe{eq:kineticupre}, $f_0(x,\xi)=\mathbf{1}_{u_0(x)>\xi}$. We have used the brackets $\<\cdot,\cdot\>$ to denote the duality between $C^\infty_c(\T^N\times\R)$ and the space of distributions over $\T^N\times\R$. In what follows, we will denote similarly the integral 
\begin{equation*}
\<F,G\>=\int_{\T^N}\int_\R F(x,\xi)G(x,\xi) dx d\xi,\quad F\in L^p(\T^N\times\R), G\in L^q(\T^N\times\R),
\end{equation*}
where $1\leq p\leq +\infty$ and $q$ is the conjugate exponent of $p$. In \refe{eq:kineticupre} also, we have indicated the dependence of $g_k$ and $\GG^2$ on $u$, which is actually absent in the additive case and we have used (with $\phi=\partial_\xi\varphi$) the shorthand $m(\phi)$ for
\begin{equation*}
m(\phi)=\int_{\T^N\times[0,T]\times\R} \phi(x,t,\xi)dm(x,t,\xi),\quad \phi\in C_b(\T^N\times[0,T]\times\R).
\end{equation*}

Note that a solution $u$ in the sense of Definition \ref{defkineticsol} is not a process in the 
usual sense since it is only defined almost everywhere with respect to the time. Part of our 
work below is to show that  $u$ has a natural representative which has almost sure continuous
trajectories with values in $L^p(\T^N)$.

Equation~\refe{eq:kineticupre} is the weak form of the equation
\begin{equation}
(\partial_t+a(\xi)\cdot\nabla)\mathbf{1}_{u>\xi}=\delta_{u=\xi}\Phi\dot W+\partial_\xi(m-\frac{1}{2}\GG^2\delta_{u=\xi}).
\label{eq:kineticD}\end{equation}
We present now a {\it formal} derivation of equation \refe{eq:kineticD} from \refe{stoSCL} {\it in the case} $m=0$ (see also Section~\ref{sec:viscousappx}, where we give a rigorous derivation of the kinetic formulation at the level of the viscous approximation): it is essentially a consequence of It\^o Formula. Indeed, by the identity 
$(\mathbf{1}_{u>\xi},\theta'):=\int_\R \mathbf{1}_{u>\xi}\theta'(\xi) d\xi=\theta(u)-\theta(-\infty)$, satisfied for $\theta\in C^\infty(\R)$,
and by It\^o Formula, we have
\begin{align*}
d(\mathbf{1}_{u>\xi},\theta')
	&=\theta'(u)(-a(u)\cdot\nabla u dt+\Phi(u)d W)+\frac{1}{2}\theta''(u)\GG^2 dt\\
	&=-\div(\int^u a(\xi)\theta'(\xi)d\xi)dt+\frac{1}{2}\theta''(u)\GG^2 dt+\theta'(u) \Phi(u)dW\\
	&=-\div((a\mathbf{1}_{u>\xi},\theta'))dt-\frac{1}{2} (\partial_\xi(\GG^2\delta_{u=\xi}),\theta')dt+(\delta_{u=\xi},\theta'\Phi dW).
\end{align*}
Taking $\theta(\xi)=\int_{-\infty}^\xi\varphi$, we then obtain the kinetic formulation with $m=0$.
The measure $m$ is sometimes (quite improperly if no action, or Lagrangian, is precisely defined) interpreted as a Lagrange multiplier for the evolution of $f$ by $\partial_t+a\cdot\nabla$ under the constraint $f=\mathrm{graph}=\mathbf{1}_{u>\xi}$. It comes into play only when $u$ becomes discontinuous (occurrence of shocks); in particular, it does not appear in the computation above that requires some regularity of $u$ with respect to $x$ to apply the chain-rule of differentiation.

\subsection{Generalized solutions}

With the purpose to prepare the proof of existence of solution, we introduce the following definitions.

\begin{definition}[Young measure] Let $(X,\lambda)$ be a finite measure space. Let $\mathcal{P}_1(\R)$ denote the set of probability measures on $\R$. We say that a map $\nu\colon X\to\mathcal{P}_1(\R)$ is a Young measure on $X$ if, for all $\phi\in C_b(\R)$, the map $z\mapsto \nu_z(\phi)$ from $X$ to $\R$ is measurable. We say that a Young measure $\nu$ vanishes at infinity if, for every $p\geq 1$, 
\begin{equation}
\int_X\int_\R |\xi|^p d\nu_z(\xi)d\lambda(z)<+\infty.
\label{nuvanish}\end{equation}
\end{definition}

\begin{definition}[Kinetic function] Let $(X,\lambda)$ be a finite measure space. A measurable function $f\colon X\times\R\to[0,1]$ is said to be a kinetic function if there exists a Young measure $\nu$ on $X$ that vanishes at infinity such that, for $\lambda$-a.e. $z\in X$, for all $\xi\in\R$,
\begin{equation*}
f(z,\xi)=\nu_{z}(\xi,+\infty).
\end{equation*}
We say that $f$ is an {\rm equilibrium} if there exists a measurable function $u\colon X\to\R$ such that $f(z,\xi)=\mathbf{1}_{z>\xi}$ a.e., or, equivalently, $\nu_z=\delta_{u(z)}$ for a.e. $z\in X$.
\end{definition}

If $f\colon X\times\R\to[0,1]$ is a kinetic function, we denote by $\bar f$ the {\it conjugate function} $\bar f:=1-f$.

We also denote by $\chi_f$ the function defined by $\chi_f(z,\xi)=f(z,\xi)-\mathbf{1}_{0>\xi}$. Contrary to $f$, this modification is integrable. Actually, it is decreasing faster than any power of $\xi$ at infinity. Indeed, 
\begin{equation*}
\chi_f(z,\xi) = \left\{
\begin{array}{l}
\ds-\int_{(-\infty,\xi]} d\nu_z, \; \xi <0,\\
\\
\ds \int_{(\xi,+\infty)} d\nu_z, \; \xi >0.
\end{array}
\right.
\end{equation*}
Therefore
\begin{equation}
\label{e10}
|\xi|^p \int_{X} |\chi_f(z,\xi)|d\lambda(z) \le \int_{X}\int_\R |\zeta|^p d\nu_{x,t}(\zeta)
d\lambda(z) <\infty,
\end{equation}
for all $\xi\in\R$, $1\leq p<+\infty$.
\medskip

We have the following compactness results (the proof is classical and reported to appendix).

\begin{theorem}[Compactness of Young measures] Let $(X,\lambda)$ be a  finite  measured 
space such that $L^1(X)$ is separable. Let $(\nu^n)$ be a sequence of Young measures on $X$ satisfying  \refe{nuvanish} uniformly for some $p\ge 1$:
\begin{equation}
\sup_n\int_X\int_\R |\xi|^p d\nu^n_z(\xi)d\lambda(z)<+\infty.
\label{nuvanishn}\end{equation}
Then there exists a Young measure $\nu$ on $X$ and a subsequence still denoted $(\nu^n)$ such that, for all $h\in L^1(X)$, for all $\phi\in C_b(\R)$,
\begin{equation}
\lim_{n\to+\infty}\int_X h(z) \int_\R \phi(\xi) d\nu^n_z(\xi) d\lambda(z)= \int_X h(z) \int_\R \phi(\xi) d\nu_z(\xi) d\lambda(z).
\label{cvYoungMeasure}\end{equation}
\label{th:youngmeasure}\end{theorem}
 
\begin{corollary}[Compactness of kinetic functions] Let $(X,\lambda)$ be a  finite  measured 
space such that $L^1(X)$ is separable. Let $(f_n)$ be a sequence of kinetic functions on $X\times\R$: $f_n(z,\xi)=\nu^n_z(\xi,+\infty)$ where $\nu^n$ are Young measures on $X$ satisfying \refe{nuvanishn}. Then there exists a kinetic function $f$ on $X\times\R$ such that $f_n\rightharpoonup f$ in $L^\infty(X\times\R)$ weak-*.
\label{cor:kineticfunctions}\end{corollary}

We will also need the following result.

\begin{lemma}[Convergence to an equilibrium] Let $(X,\lambda)$ be a finite measure space. Let $p> 1$. Let $(f_n)$ be a sequence of kinetic functions on $X\times\R$: $f_n(z,\xi)=\nu^n_z(\xi,+\infty)$ where $\nu^n$ are Young measures on $X$ satisfying \refe{nuvanishn}. Let $f$ be a kinetic function on $X\times\R$ such that $f_n\rightharpoonup f$ in $L^\infty(X\times\R)$ weak-*. Assume that $f_n$ and $f$ are equilibria: 
$$
f_n(z,\xi)=\mathbf{1}_{u_n(z)>\xi},\quad f(z,\xi)=\mathbf{1}_{u(z)>\xi}.
$$
Then, for all $1\leq q<p$, $u_n\to u$ in $L^q(X)$ strong.
\label{lem:weakstrongeq}\end{lemma}

Note that if $f$ is a kinetic function then $\partial_\xi f=-\nu$ is non-negative.
Observe also that, in the context of Definition~\ref{defkineticsol}, setting $f=\mathbf{1}_{u>\xi}$, we have $\partial_\xi f=-\delta_{u=\xi}$ and $\nu:=\delta_{u=\xi}$ is a Young measure on $\Omega\times\T^N\times(0,T)$. The measure $\nu$ vanishes at infinity (it even satisfies the stronger condition \refe{eq:integrabilityf} below). Therefore any solution will also be a {\it generalized solution}, according to the definition below.

\begin{definition}[Generalized solution]
\label{d4} Let $f_0\colon\Omega\times\T^N\times\R\to[0,1]$ be a kinetic function. A measurable function $f\colon\Omega\times\T^N\times[0,T]\times\R\to[0,1]$ is said to be a generalized solution to~\refe{stoSCL} with initial datum $f_0$ if $(f(t))$ is predictable and is a kinetic function such that: for all $p\geq 1$, $\nu:=-\partial_\xi f$ satisfies 
\begin{equation}
\E\left(\esssup_{t\in[0,T]}\int_{\T^N}\int_\R|\xi|^p d\nu_{x,t}(\xi) dx\right) \leq C_p,
\label{eq:integrabilityf}\end{equation}
where $C_p$ is a positive constant and: there exists a kinetic measure $m$ such that for all $\varphi\in C^1_c(\T^N\times[0,T)\times\R)$, 
\begin{equation}
\label{eq:kineticfpre}
\begin{array}{l}
\ds \int_0^T\<f(t),\partial_t \varphi(t)\>dt+\<f_0,\varphi(0)\>
+\int_0^T \<f(t),a(\xi)\cdot\nabla\varphi(t)\>dt\\
\ds =-\sum_{k\geq 1}\int_0^T\int_{\T^N}\int_\R g_k(x,\xi)\varphi(x,t,\xi)d\nu_{x,t}(\xi)dxd\beta_k(t)\\
\ds -\frac{1}{2}\int_0^T\int_{\T^N}\int_\R \partial_\xi\varphi(x,t,\xi)\GG^2(x,\xi)d\nu_{(x,t)}(\xi) dx dt+m(\partial_\xi\varphi), \mbox{ a.s.}
\end{array}
\end{equation}
\end{definition}

Observe that, if $f$ is a generalized solution such that $f=\mathbf{1}_{u>\xi}$, then 
$u(t,x)=\int_\R \chi_f(x,t,\xi)d\xi$, hence $u$ is predictable. Moreover,
$\nu=\delta_{u=\xi}$ and
$$
\int_{\T^N} |u(t,x)|^p dx = \int_{\T^N}\int_\R |\xi|^p d\nu_{x,t}(\xi)dx .
$$
Condition~\refe{eq:integrabilityu} is thus contained in the condition~\refe{eq:integrabilityf}. 
\medskip

We conclude this paragraph with two remarks. The first remark is the following
\begin{lemma}[Distance to equilibrium] Let $(X,\lambda)$ be a finite measure space. Let $f\colon X\times\R\to[0,1]$ be a kinetic function. Then
$$
m(\xi):=\int_{-\infty}^\xi (\mathbf{1}_{u>\zeta}-f(\zeta))d\zeta,\quad\mbox{where } u:=\int_\R\chi_{f}(\zeta)d\zeta,
$$
is well defined and non-negative.
\label{lemdisteq}\end{lemma}

Note in particular that the difference $f(\xi)-\mathbf{1}_{u>\xi}$ writes $\partial_\xi m$ where $m\geq 0$.
\medskip

{\bf Proof of Lemma~\ref{lemdisteq}: } Let $\nu_z=-\partial_\xi f(z,\cdot)$, $z\in X$. By Jensen's inequality, we have
\begin{equation}\label{JensenH}
H\left(\int_\R \zeta d\nu_z(\zeta)\right)\leq \int_\R H(\zeta)d\nu_z(\zeta)
\end{equation}
for all convex sub-linear function $H\colon\R\to\R$. Note that
$$
u(z)=\int_\R f(z,\zeta)-\mathbf{1}_{0>\zeta} d\zeta=\int_\R \zeta d\nu_z(\zeta)
$$
by integration by parts. By integration by parts, we also have, for $H\in C^1(\R)$ and sub-linear,
$$
\int_\R H(\zeta)d\nu_z(\zeta)=H(0)+\int_\R H'(\zeta)(f(z,\zeta)-\mathbf{1}_{0>\zeta} )d\zeta
$$
and
$$
H(u(z))=\int_\R H(\zeta)d\delta_{u(z)}(\zeta)=H(0)+\int_\R H'(\zeta)(\mathbf{1}_{u(z)>\zeta}-\mathbf{1}_{0>\zeta} )d\zeta.
$$
By \refe{JensenH}, it follows that
$$
\int_\R H'(\zeta)(f(z,\zeta)-\mathbf{1}_{u(z)>\zeta} )d\zeta\geq 0
$$
for all convex and sub-linear $H\in C^1(\R)$. Approximating $\zeta\mapsto (\zeta-\xi)^-$ by such functions $H$, we obtain $m(\xi)\geq 0$. \qed
\bigskip

Our second remark is about the time continuity of the solution  (see also \cite{CancesGallouet10} and references therein on this subject). Generalized solutions are a useful and natural tool for the analysis of {\it weak} solutions to \refe{stoSCL}, {\it i.e.} solutions that are weak with respect to space {\it and} time, but the process of relaxation that generalizes the notion of solution introduces additional difficulties regarding the question of time continuity of solutions. To illustrate this fact, let us consider for example the following equation (the ``Collapse" equation in the Transport-Collapse method of Brenier~\cite{Brenier81,Brenier83})
\begin{equation}
\partial_t f(t)=\mathbf{1}_{u(t)>\xi}-f,\quad u(t):=\int_\R\chi_{f(t)}(\xi)d\xi,
\label{Collapse}\end{equation}
with initial datum $f_0(\xi)$ a kinetic function. Integrating \refe{Collapse} with respect to $\xi$ shows that $u=u_0$ is constant and gives
\begin{equation*}
f(t)=e^{-t}f_0+(1-e^{-t})\mathbf{1}_{u_0>\xi},
\end{equation*}
{\it i.e.} $f(t)$ is describing the progressive and continuous ``collapse" from $f_0$ to $\mathbf{1}_{u_0>\xi}$. By Lemma~\ref{lemdisteq}, 
\begin{equation*}
m(t,\xi):=\int_{-\infty}^\xi (\mathbf{1}_{u>\zeta}-f(t,\zeta)) d\zeta \geq 0
\end{equation*}
for all $t,\xi$. More generally,
\begin{equation}
\int_{-\infty}^\xi (f(\tau,\zeta)-f(t,\zeta)) d\zeta \geq 0,\quad\forall \tau>t,\forall\xi\in\R,
\label{posmtaut}\end{equation}
as we obtain by integration of \refe{Collapse} with respect to $s\in(t,\tau)$ and $\zeta<\xi$. Observe also that $f$ satisfies $\partial_t f=\partial_\xi m$, $m\geq 0$.  Now erase an interval $[t_1,t_2]$ in the evolution of $f$. Then
\begin{equation*}
g(t)=\hat f(t):=f(t)\mathbf{1}_{[0,t_1]}(t)+f(t+t_2-t_1)\mathbf{1}_{(t_1,+\infty)}(t)
\end{equation*}
satisfies 
\begin{align*}
\partial_t g&=\partial_\xi \hat m+(f(t_2)-f(t_1))\delta(t-t_1)\\
&=\partial_\xi n,\quad n(t,\xi):=\hat m(t,\xi)+\int_{-\infty}^\xi (f(t_2,\zeta)-f(t_1,\zeta))d\zeta\delta(t-t_1).
\end{align*}
By \refe{posmtaut}, $n$ is non-negative, but, unless $f_0=\mathbf{1}_{u_0>\xi}$ in which case $f$ is
constant and $g=f$, $g$ is discontinuous at $t=t_1$. In the analysis of a generalized solution $f$, we thus first show the existence of modifications $f^+$ and $f^-$ of $f$ being respectively right- and left-continuous everywhere and we work on $f^\pm$ in most of the proof of uniqueness and reduction (Theorem~\ref{th:Uadd}). Finally, we obtain the time-continuity of solutions in Corollary~\ref{cor:timecontinuity}.

\subsection{Left and right limits of generalized solution}

We show in the following proposition that, almost surely, any generalized solution admits possibly different left and right weak limits at any point $t\in[0,T]$. This property is important to prove a comparison principle which allows to prove uniqueness. Also, it allows us to see that the weak 
form \eqref{eq:kineticfpre} of the equation satisfied by  
a generalized solution can be strengthened. We write below a  formulation which is weak only with respect to $x$
and $\xi$.

Note that we obtain continuity with respect to time of solutions in Corollary~\ref{cor:timecontinuity} below. 

\begin{proposition}[Left and right weak limits] Let $f_0$ be a kinetic initial datum. Let $f$ be a generalized solution to~\refe{stoSCL} with initial datum $f_0$. Then  $f$ admits almost surely left and right limits at all point $t_*\in [0,T]$. More precisely, for all $t_*\in [0,T]$ there exists some kinetic functions $f^{*,\pm}$ on $\Omega\times\T^N\times\R$ such that $\P$-a.s. 
\begin{equation*}
\<f(t_*-\eps),\varphi\> \to \<f^{*,-},\varphi\>
\end{equation*}
and
\begin{equation*}
\<f(t_*+\eps),\varphi\> \to \<f^{*,+},\varphi\>
\end{equation*}
as $\eps\to0$ for all $\varphi\in C^1_c(\T^N\times\R)$. Moreover, almost surely, 
\begin{equation}
\<f^{*,+}-f^{*,-},\varphi\>=-\int_{\T^N\times[0,T]\times\R}\partial_\xi\varphi(x,\xi)\mathbf{1}_{\{t_*\}}(t) dm(x,t,\xi).
\label{f+VSf-}\end{equation}
In particular, almost surely, the set of $t_*\in[0,T]$ such that $f^{*,-}\ne f^{*,+}$ is countable.
\label{prop:LRlimits}\end{proposition}

In the following, for a  generalized solution $f$, we define $f^\pm$ by $f^\pm(t_*)=f^{*\pm}$, 
$t_*\in [0,T]$. Note that, since we are dealing with a filtration associated to brownian motions, 
$f^\pm$ are also predictable. Also $f=f^+=f^-$ almost everywhere in time and 
we can take any of them in an integral 
with respect to the Lebesgue measure or in a stochastic integral.
On the contrary, if the integration is with respect to a measure - typically a kinetic measure in this article -,
the integral is not well defined for $f$ and may differ if one chooses $f^+$ or $f^-$.

\smallskip

{\bf Proof of Proposition~\ref{prop:LRlimits}.}Without loss of generality, we assume that
$\Omega=C([0,T];U)$, where $U$ is a Hilbert space and $H\subset U$ with Hilbert-Schmidt 
embedding, $\F$ is the Borel $\sigma$-algebra of $C([0,T],U)$ and that $\P$ is the Wiener measure on $\Omega$. 

The set of test functions $C^1_c(\T^N\times\R)$ (endowed with the topology of the uniform 
convergence on any compact of the functions and their first derivatives) is separable and we fix 
a dense countable subset $\mathcal D_1$. For all $\varphi\in C^1_c(\T^N\times\R)$, a.s., the 
map
\begin{multline}
J_\varphi\colon t\mapsto\int_0^t \<f(s),a(\xi)\cdot\nabla\varphi\>ds\\-\sum_{k\geq 1}\int_0^t\int_{\T^N}\int_\R g_k(x,\xi)\varphi(x,\xi)d\nu_{x,s}(\xi)dxd\beta_k(s)\\
+\frac{1}{2}\int_0^t\int_{\T^N}\int_\R \partial_\xi\varphi(x,\xi)\GG^2(x,\xi)d\nu_{x,s}(\xi) dx ds
\label{def:Jphi}\end{multline}
is continuous on $[0,T]$. Consequently: a.s., for all $\varphi\in \mathcal D_1$, $J_\varphi$ is continuous on $[0,T]$. 

For test functions of the form $(x,t,\xi)\mapsto\varphi(x,\xi)\alpha(t)$, $\alpha\in C^1_c([0,T])$, $\varphi\in \mathcal D_1$, Fubini Theorem and the weak formulation~\refe{eq:kineticfpre} give
\begin{equation}
\int_0^T g_\varphi(t)\alpha'(t)dt+\<f_0,\varphi\>\alpha(0)=\<m,\partial_\xi\varphi\>(\alpha),
\label{eq:weak0}\end{equation}
where $g_\varphi(t):=\<f(t),\varphi\>-J_\varphi(t)$.
This shows that $\partial_t g_\varphi$ is a Radon measure on $(0,T)$, {\it i.e.} the function $g_\varphi\in BV(0,T)$. In particular it admits left and right limits at all points $t_*\in[0,T]$. Since $J_\varphi$ is continuous, this also holds for $\<f,\varphi\>$: for all $t_*\in[0,T]$, the limits
\begin{equation*}
\<f,\varphi\>(t_*+):=\lim_{t\downarrow t_*}\<f,\varphi\>(t)\mbox{ and }\<f,\varphi\>(t_*-):=\lim_{t\uparrow t_*}\<f,\varphi\>(t)
\end{equation*}
exist. Note that:
\begin{equation*}
\<f,\varphi\>(t_*+)=\lim_{\eps\to 0}\frac{1}{\eps}\int_{t_*}^{t_*+\eps}\<f,\varphi\>(t)dt,\quad  
\<f,\varphi\>(t_*-)=\lim_{\eps\to 0}\frac{1}{\eps}\int_{t_*-\eps}^{t_*} \<f,\varphi\>(t)dt.
\end{equation*}

Let $(\eps_n)\downarrow 0$. Set $X=\Omega\times\T^N\times\R$ and let $\lambda$ denote the product measure of the Wiener measure $\P$ and of the Lebesgue measure on $\T^N\times\R$. The function 
$$
f_n:=\frac{1}{\eps_n}\int_{t_*}^{t_*+\eps_n}f(t)dt
$$
is a kinetic function which, thanks to~\refe{eq:integrabilityf}, satisfies the condition \refe{nuvanishn}. Clearly, the Borel $\sigma$ field on $\Omega\times \T^n\times\R
=C([0,T]:U)\times \T^N \times \R$ is countably generated and by \cite{cohn}, Proposition 3.4.5, $L^1(\Omega\times \T^n\times\R)$ is separable.  By Corollary~\ref{cor:kineticfunctions}, there exist a kinetic functions $f^{*,\pm}$ on $\Omega\times\T^N\times\R$ and subsequences $(\eps_{n_k^\pm})$ such that
\begin{equation*}
\frac{1}{\eps_{n_k^-}}\int_{t_*-\eps_{n_k^-}}^{t_*}f(t)dt\rightharpoonup f^{*,-},\quad \frac{1}{\eps_{n_k^+}}\int_{t_*}^{t_*+\eps_{n_k^+}}f(t)dt\rightharpoonup f^{*,+}
\end{equation*}
weakly-$*$ in $L^\infty(\Omega\times\T^N\times\R)$ as $k\to+\infty$. We deduce:
\begin{equation*}
\<f,\varphi\>(t_*+)=\<f^{*,+},\varphi\> \mbox{ and }\<f,\varphi\>(t_*-)= \<f^{*,-},\varphi\>.
\end{equation*}
Taking for $\alpha$ the hat function $\alpha(t)=\ds\frac{1}{\eps}\min((t-t_*+\eps)^+,(t-t_*-\eps)^-)$ in \refe{eq:weak0}, we obtain \refe{f+VSf-} at the limit $[\eps\to0]$. In particular, almost surely, $f^{*,+}=f^{*,-}$ whenever $m$ has no atom at $t_*$.

We thus have proved the result for $\varphi\in\mathcal D_1$. Since $\mathcal D_1$ is dense in 
$C^1_c(\T^N\times\R)$, it is easy to see that in fact everything holds a.s. for every $\varphi\in C^1_c(\T^N\times\R)$. \qed

\begin{remark}[Uniform bound] Note that, by construction, $f^\pm$ satisfy the bound \refe{eq:integrabilityf} uniformly in time:
\begin{equation}
\E\left(\sup_{t\in[0,T]}\int_{\T^N}\int_\R|\xi|^p d\nu^\pm_{x,t}(\xi) dx\right) \leq C_p,
\label{eq:integrabilityfpm}\end{equation}
\label{rk:integrabilityfpm}\end{remark}

Once we have proved the existence of left and right limits everywhere, we can derive a kinetic formulation at given $t$ ({\it i.e.} weak in $(x,\xi)$ only). Taking in \refe{eq:kineticfpre} a test function of the form $(x,s,\xi)\mapsto\varphi(x,\xi)\alpha(s)$ where $\alpha$ is the function
$$
\alpha(s)=\left\{
\begin{array}{ll}
1,\quad &s\le t,\\
\ds1-\frac{s-t}\eps, \quad &t\le s\le t+\eps,\\
0,& t+\eps \le s,
\end{array}
\right.
$$
we obtain  at the limit $[\eps\to0]$: for all $t\in [0,T]$ and $\varphi\in C^1_c(\T^N\times\R)$,
\begin{equation}
\begin{array}{l}
\ds -\<f^+(t),\varphi\>+\<f_0,\varphi\>
+\int_0^t \<f(s),a(\xi)\cdot\nabla\varphi\>ds\\
\ds =-\sum_{k\geq 1}\int_0^t\int_{\T^N}\int_\R g_k(x,\xi)\varphi(x,\xi)d\nu_{x,s}(\xi)dxd\beta_k(s)\\
\ds -\frac{1}{2}\int_0^t\int_{\T^N}\int_\R \partial_\xi\varphi(x,\xi)\GG^2(x,\xi)d\nu_{(x,s)}(\xi) dx ds+
\<m,\partial_\xi\varphi\>([0,t]), \mbox{ a.s.},
\end{array}
\label{eq:19}\end{equation}
where 
$\ds \<m,\partial_\xi\varphi\>([0,t])= \int_{\T^N\times[0,t]\times\R} \partial_\xi\varphi(x,\xi)dm(x,s,\xi)$.
\medskip

\begin{remark}[The case of equilibrium] Assume that $f^{*,-}$ is at equilibrium in \refe{f+VSf-}: there exists a random variable $u^*\in L^1(\T^d)$ such that $f^{*,-}=\mathbf{1}_{u^*>\xi}$ a.s. Let $m^*$ denote the restriction of $m$ to $\T^N\times\{t_*\}\times\R$. We thus have
\begin{equation*}
f^{*,+}-\mathbf{1}_{u^*>\xi}=\partial_\xi m^*.
\end{equation*}
In particular, by the condition at infinity \refe{inftym} on $m$ the integral of the rhs vanishes and we have: almost surely, for a.e. $x\in\T^N$, 
$$
\int_\R( {f^{*,+}}(x,\xi)-\mathbf{1}_{0>\xi}) d\xi=\int_\R(\mathbf{1}_{u^*>\xi}-\mathbf{1}_{0>\xi}) d\xi=u^*.
$$
By Lemma~\ref{lemdisteq},
$$
p^*\colon\xi\mapsto\int_{-\infty}^\xi (\mathbf{1}_{u^*>\zeta}-f^{*,+}(\zeta))d\zeta
$$
is non-negative. Besides, $\partial_\xi (m^*+p^*)=0$, hence $m^*+p^*$ is constant, and actually vanishes by the condition at infinity \refe{inftym} and the obvious fact that $p$ also vanishes when $|\xi|\to+\infty$. Since $m^*,p^*\geq 0$, we finally obtain $m^*=0$ and $f^{*,+}=f^{*,-}$: in conclusion, when $f^{*,-}$ is at equilibrium, \refe{f+VSf-} is trivial and we have no discontinuity at $t_*$.
\label{RKequilibrium}\end{remark}
\section{Comparison, uniqueness, entropy solution and regularity}\label{sec:comparison}

\subsection{Doubling of variables}
In this paragraph, we prove a technical proposition relating two generalized solutions $f_i$, $i=1,2$ of the equation
\begin{equation}\label{stoSCLi}
du_i+\div(A(u_i))dt=\Phi(u_i)dW.
\end{equation}

\begin{proposition} Let $f_i$, $i=1,2$, be generalized solution to \refe{stoSCLi}. Then, for $0\leq t\leq T$, and non-negative test functions $\rho\in C^\infty(\T^N)$, $\psi\in C^\infty_c(\R)$, we have 
\begin{multline}
\E\int_{(\T^N)^2}\int_{\R^2} \rho(x-y)\psi(\xi-\zeta)f_1^\pm(x,t,\xi)\bar f_2^\pm(y,t,\zeta) d\xi d\zeta dx dy \\
\leq\E\int_{(\T^N)^2}\int_{\R^2} \rho(x-y)\psi(\xi-\zeta) f_{1,0}(x,\xi)\bar f_{2,0}(y,\zeta)d\xi d\zeta dx dy +\mathrm{I}_\rho+\mathrm{I}_\psi,
\label{CR0}\end{multline}
where
\begin{multline*}
\mathrm{I}_\rho=\E\int_0^t\int_{(\T^N)^2}\int_{\R^2} f_1(x,s,\xi)\bar f_2(y,s,\zeta)(a(\xi)-a(\zeta))\psi(\xi-\zeta)d\xi d\zeta\\
 \cdot\nabla_x\rho(x-y) dx dy ds
\end{multline*}
and
\begin{multline*}
\mathrm{I}_\psi=\frac{1}{2}\int_{(\T^N)^2}\rho(x-y)\E\int_0^t\int_{\R^2} \psi(\xi-\zeta)\\
\times\sum_{k\geq 1}|g_{k}(x,\xi)-g_{k}(y,\zeta)|^2 d\nu^1_{x,s}\otimes\nu^2_{y,s}(\xi,\zeta) dx dy ds.
\end{multline*}
\label{prop:CR}\end{proposition}

\begin{remark}\label{intfbarf} Each term in \refe{CR0} is finite. Let us for instance consider the first one on the right hand side. Let us introduce the auxiliary functions 
\begin{equation*}
\psi_1(\xi)=\int_{-\infty}^\xi\psi(s) ds,\quad \psi_2(\zeta)=\int_{-\infty}^\zeta\psi_1(\xi)d\xi,
\end{equation*}
which are well-defined  since $\psi$ is compactly supported. Note that both $\psi_1$ and $\psi_2$ vanish at $-\infty$. When $\xi\to+\infty$, $\psi_1$ remains bounded while $\psi_2$ has linear 
growth. 
To lighten notations, we omit the index $0$. Let us set $\bar f_2=1-f_2$. 
In the case where $f_1$ and $f_2$ correspond to kinetic solutions, {\it i.e.}
$f_i=\mathbf{1}_{u_i>\xi}$, we compute (forgetting the dependence upon $t$ and $x$): $\bar f_2(\zeta)=\mathbf{1}_{u_2\leq\zeta}$ and
\begin{equation*}
\int_{\R^2}\psi(\xi-\zeta)f_1(\xi)\bar f_2(\zeta)d\xi d\zeta=\psi_2(u_1-u_2).
\end{equation*}
In the case of generalized solutions, we introduce the integrable modifications $\chi_{f_i}$ of $f_i$, $i=1,2$: 
\begin{equation*}
f_1(\xi)=\chi_{f_1}(\xi)+\mathbf{1}_{0>\xi},\quad \bar f_2(\zeta)=\mathbf{1}_{0\leq\zeta}-\chi_{f_2}(\zeta). 
\end{equation*}
Accordingly, we have, by explicit integration:
\begin{multline}
\int_{\R^2}\psi(\xi-\zeta)f_1(\xi)\bar f_2(\zeta)d\xi d\zeta=-\int_{\R^2}\psi(\xi-\zeta)\chi_{f_1}(\xi)\chi_{f_2}(\zeta)d\xi d\zeta\\
+\int_{\R}\psi_1(\xi)\chi_{f_1}(\xi)d\xi
-\int_{\R}\psi_1(\zeta)\chi_{f_2}(-\zeta)d\zeta+\psi_2(0)
\label{remarkfinite}\end{multline}
Each term in the right hand-side of \refe{remarkfinite} is indeed finite by \refe{e10}.
\end{remark}

{\bf Proof of Proposition~\ref{prop:CR}:} Set $G_1^2(x,\xi)=\sum_{k=1}^\infty|g_{k,1}(x,\xi)|^2$ and $G_2^2(y,\zeta)=\sum_{k=1}^\infty|g_{k,2}(y,\zeta)|^2$. Let $\varphi_1\in C^\infty_c(\T^N_x\times\R_\xi)$ and $\varphi_2\in C^\infty_c(\T^N_y\times\R_\zeta)$. By \eqref{eq:19}, we have
$$
 \<f_1^{+}(t),\varphi_1\>=\<m_1^*,\partial_\xi\varphi_1\>([0,t])+F_1(t)
$$
with
$$
F_1(t)=\sum_{k\geq 1}\int_0^t\int_{\T^N}\int_\R g_{k,1}\varphi_1 d\nu^1_{x,s}(\xi)dxd\beta_k(s)
$$
and
\begin{multline*}
\<m_1^*,\partial_\xi\varphi_1\>([0,t])= \<f_{1,0},\varphi_1\>\delta_0([0,t])+
\int_0^t \<f_1,a\cdot\nabla\varphi_1\>ds\\
+\frac{1}{2}\int_0^t\int_{\T^N}\int_\R \partial_\xi\varphi_1\GG^2_1 d\nu^1_{(x,s)}(\xi) dx ds-\<m_1,\partial_\xi\varphi_1\>([0,t]).
\end{multline*}
Note that, by Remark~\ref{RKequilibrium}, $\<m_1,\partial_\xi\varphi_1\>(\{0\})=0$ and thus the value of $\<m_1^*,\partial_\xi\varphi_1\>(\{0\})$ is $ \<f_{1,0},\varphi_1\>$. Similarly
$$
 \<\bar f_2^{+}(t),\varphi_2\>= \<\bar m_2^*,\partial_\zeta\varphi_2\>([0,t])+\bar F_2(t)
$$
with
$$
\bar F_2(t)=-\sum_{k\geq 1}\int_0^t\int_{\T^N}\int_\R g_{k,2}\varphi_2 d\nu^2_{y,s}(\zeta)dyd\beta_k(s)
$$
and
\begin{multline*}
\<\bar m_2^*,\partial_\zeta\varphi_2\>([0,t])=\<\bar f_{2,0},\varphi_2\>\delta_0([0,t]) +
\int_0^t \<\bar f_2,a\cdot\nabla\varphi_2\>ds\\
-\frac{1}{2}\int_0^t\int_{\T^N}\int_\R \partial_\xi\varphi_2\GG^2_1 d\nu^2_{(y,s)}(\zeta) dy ds
+\<m_2,\partial_\zeta\varphi_2\>([0,t]),
\end{multline*}
where $\<\bar m_2^*,\partial_\zeta\varphi_2\>(\{0\})=\<\bar f_{2,0},\varphi_2\>$. Let $\alpha(x,\xi,y,\zeta)=\varphi_1(x,\xi)\varphi_2(y,\zeta)$. 
Using It\^o formula for $F_1(t)\bar F_2(t)$, integration by parts for functions of finite variation (see for instance \cite{RevuzYor99}, chapter 0) for $\<m_1^*,\partial_\xi\varphi_1\>([0,t])\<\bar m_2^*,\partial_\zeta\varphi_2\>([0,t])$, which gives
\begin{multline*}
\<m_1^*,\partial_\xi\varphi_1\>([0,t])\<\bar m_2^*,\partial_\zeta\varphi_2\>([0,t])\\
=\<m_1^*,\partial_\xi\varphi_1\>(\{0\})\<\bar m_2^*,\partial_\zeta\varphi_2\>(\{0\})
+\int_{(0,t]}\<m_1^*,\partial_\xi\varphi_1\>([0,s))d\<\bar m_2^*,\partial_\zeta\varphi_2\>(s)\\
+\int_{(0,t]}\<\bar m_2^*,\partial_\zeta\varphi_2\>([0,s])d\<m_1^*,\partial_\xi\varphi_1\>(s)
\end{multline*}
and the following formula
$$
\<m_1^*,\partial_\xi\varphi_1\>([0,t])\bar F_2(t)
= \int_0^t \<m_1^*,\partial_\xi\varphi_1\>([0,s]) d\bar F_2(s)
+\int_0^t \bar F_2(s) \<m_1^*,\partial_\xi\varphi_1\>(ds),
$$
which is easy to obtain since $\bar F_2$ is continuous, and a similar formula 
for $\<\bar m_2^*,\partial_\zeta\varphi_2\>([0,t])\bar F_1(t)$, we obtain
that 
\begin{equation*}
\<f_1^+(t),\varphi_1\>\<\bar f_2^+(t),\varphi_2\>=\<\<f_1^+(t)\bar f_2^+(t),\alpha\>\>
\end{equation*}
satisfies
\begin{multline}
\E \<\<f_1^+(t)\bar f_2^+(t),\alpha\>\>
=\<\<f_{1,0}\bar f_{2,0},\alpha\>\>\\
+\E\int_0^t\int_{(\T^N)^2}\int_{\R^2} f_1 \bar f_2 (a(\xi)\cdot\nabla_x +a(\zeta)\cdot\nabla_y)\alpha d\xi d\zeta dx dy ds\\
+\frac12 \E\int_0^t\int_{(\T^N)^2}\int_{\R^2} \partial_\xi\alpha \bar f_2(s)\GG^2_1 d\nu^1_{(x,s)}(\xi) d\zeta dx dy ds\\
-\frac12 \E\int_0^t\int_{(\T^N)^2}\int_{\R^2}\int_\R \partial_\zeta\alpha f_1(s) \GG^2_2 d\nu^2_{(y,s)}(\zeta) d\xi dy dx ds\\
-\E\int_0^t\int_{(\T^N)^2}\int_{\R^2}\GG_{1,2} \alpha d\nu^1_{x,s}(\xi) d\nu^2_{y,s}(\zeta) dx dy\\
-\E\int_{(0,t]}\int_{(\T^N)^2}\int_{\R^2} \bar f_2^+(s)\partial_\xi\alpha dm_1(x,s,\xi) d\zeta dy\\
+\E\int_{(0,t]}\int_{(\T^N)^2}\int_{\R^2} f_1^-(s)\partial_\zeta\alpha dm_2(y,s,\zeta) d\xi dx
\label{CR1}\end{multline}
where $\GG_{1,2}(x,y;\xi,\zeta):=\sum_{k\geq 1}g_{k,1}(x,\xi) g_{k,2}(y,\zeta)$ and $\<\<\cdot,\cdot\>\>$ denotes the duality distribution over $\T^N_x\times\R_\xi\times\T^N_y\times\R_\zeta$. By a 
density argument, \refe{CR1} remains true for any test-function $\alpha\in C^\infty_c(\T^N_x\times\R_\xi\times\T^N_y\times\R_\zeta)$. Using similar arguments as in Remark \ref{intfbarf}, the assumption that $\alpha$ is compactly supported can be relaxed thanks to the condition at infinity \refe{inftym} on $m_i$ and \refe{nuvanish} on $\nu^i$, $i=1,2$. Using truncates of $\alpha$, we obtain that \refe{CR1} remains true if $\alpha\in C^\infty_b(\T^N_x\times\R_\xi\times\T^N_y\times\R_\zeta)$ is compactly supported in a neighbourhood of the diagonal
\begin{equation*}
\{(x,\xi,x,\xi);x\in\T^N,\xi\in\R\}.
\end{equation*} 
We then take $\alpha=\rho\psi$ where $\rho=\rho(x-y)$, $\psi=\psi(\xi-\zeta)$. Note the remarkable identities
\begin{equation}
(\nabla_x+\nabla_y)\alpha=0,\quad (\partial_\xi+\partial_\zeta)\alpha=0.
\label{ddvar}\end{equation}
In particular, the last term in \refe{CR1} is
\begin{align*}
&\E\int_{(0,t]}\int_{(\T^N)^2}\int_{\R^2} f_1^-(s)\partial_\zeta\alpha d\xi dx dm_2(y,s,\zeta)\\
=&-\E\int_{(0,t]}\int_{(\T^N)^2}\int_{\R^2} f_1^-(s)\partial_\xi\alpha d\xi dx dm_2(y,s,\zeta)\\
=&- \E\int_{(0,t]}\int_{(\T^N)^2}\int_{\R^2}\alpha  d\nu^{1,-}_{x,s}(\xi) dx dm_2(y,s,\zeta)\leq 0 
\end{align*}
since $\alpha\geq 0$. The symmetric term 
\begin{align*}
&-\E\int_{(0,t]}\int_{(\T^N)^2}\int_{\R^2} \bar f_2^+(s)\partial_\xi\alpha dm_1(x,s,\xi) d\zeta dy\\
=&-\E\int_{(0,t]}\int_{(\T^N)^2}\int_{\R^2} \alpha d\nu^{2,+}_{y,s}(\zeta) dy dm_1(x,s,\xi)
\end{align*}
is, similarly, non-positive. Consequently, we have
\begin{equation}
\E \<\<f_1^+(t)\bar f_2^+(t),\alpha\>\>\leq \<\<f_{1,0}\bar f_{2,0},\alpha\>\>+\mathrm{I}_\rho+\mathrm{I}_\psi,
\label{CR2}\end{equation}
where
\begin{equation*}
\mathrm{I}_\rho:=\E\int_0^t\int_{(\T^N)^2}\int_{\R^2} f_1 \bar f_2 (a(\xi)\cdot\nabla_x +a(\zeta)\cdot\nabla_y)\alpha d\xi d\zeta dx  dy ds
\end{equation*}
and
\begin{multline*}
\mathrm{I}_\psi=\frac12\E\int_0^t\int_{(\T^N)^2}\int_{\R^2} \partial_\xi\alpha \bar f_2(s)\GG^2_1 d\nu^1_{(x,s)}(\xi)  d\zeta dx dy ds\\
-\frac12\E\int_0^t\int_{(\T^N)^2}\int_{\R^2}\int_\R \partial_\zeta\alpha f_1(s) \GG^2_2 d\nu^2_{(x,s)}(\zeta)  d\xi dy dx ds\\
-\E\int_0^t\int_{(\T^N)^2}\int_{\R^2}\int_\R \GG_{1,2} \alpha d\nu^1_{x,s}(\xi) d\nu^2_{y,s}(\zeta) dx dy.
\end{multline*}
Equation~\refe{CR2} is indeed equation~\refe{CR0} for $f_i^+$ since, by \refe{ddvar},
\begin{equation*}
\mathrm{I}_\rho=\E\int_0^t\int_{(\T^N)^2}\int_{\R^2} f_1\bar f_2 (a(\xi)-a(\zeta))\cdot\nabla_x\alpha d\xi  d\zeta dx dy ds
\end{equation*}
and, by \refe{ddvar} also and integration by parts,
\begin{align*}
\mathrm{I}_\psi&=\frac12\E\int_0^t\int_{(\T^N)^2}\int_{\R^2} \alpha (\GG^2_1+\GG^2_2-2\GG_{1,2}) d\nu^1_{x,s}\otimes\nu^2_{y,s}(\xi,\zeta) dx dy ds\\
&=\frac12\E\int_0^t\int_{(\T^N)^2}\int_{\R^2} \alpha \sum_{k\geq 0}|g_{k}(x,\xi)-g_{k}(y,\zeta)|^2 d\nu^1_{x,s}\otimes\nu^2_{y,s}(\xi,\zeta) dx dy ds.
\end{align*}
To obtain the result for $f_i^-$, we take $t_n \uparrow t$, write \refe{CR0} for $f_i^+(t_n)$ and let
$n\to\infty$.

\subsection{Uniqueness, reduction of generalized solution}

In this section we use Proposition~\ref{prop:CR} above to deduce the uniqueness of solutions and the reduction of generalized solutions to solutions.

\begin{theorem}[Uniqueness, Reduction] Let $u_0\in L^\infty(\T^N)$. Assume~\refe{D0}-\refe{D1}. 
Then, there is at most one solution with initial datum $u_0$ to \refe{stoSCL}. Besides, any generalized 
solution $f$ is actually a solution, {\it i.e.} if $f$ is a ge\-ne\-ra\-lized solution to \refe{stoSCL} with initial 
datum $\mathbf{1}_{u_0>\xi}$, then there exists a solution $u$ to \refe{stoSCL} with initial datum 
$u_0$ such that $f(x,t,\xi)=\mathbf{1}_{u(x,t)>\xi}$ a.s., for a.e. $(x,t,\xi)$. 
\label{th:Uadd}\end{theorem}

\begin{corollary}[Continuity in time]Let $u_0\in L^\infty(\T^N)$. Assume~\refe{D0}-\refe{D1}. Then, 
for every $p\in[1,+\infty)$, the solution $u$ to \refe{stoSCL} with initial datum $u_0$ has a representative in $L^p(\Omega;L^\infty(0,T;L^p(\T^N)))$ with almost sure
continuous trajectories in $L^p(\T^N)$.
\label{cor:timecontinuity}\end{corollary}
{\bf Proof of Theorem~\ref{th:Uadd}:} Consider first the additive case: $\Phi(u)$ independent on $u$. Let $f_i
$, $i=1,2$ be two generalized solutions to \refe{stoSCL}. Then, we use \refe{CR0} with $g_k$ 
independent on $\xi$ and $\zeta$. By \refe{D1}, the last term $\mathrm{I}_\psi$ is bounded by
\begin{equation*}
\frac{t D_1}{2}\|\psi\|_{L^\infty}\int_{(\T^N)^2}|x-y|^2\rho(x-y)dx dy.
\end{equation*}
We then take $\psi:=\psi_{\delta}$ and $\rho=\rho_\eps$ where $(\psi_\delta)$ and $(\rho_\eps)$ are 
approximations to the identity on $\R$ and $\T^N$ respectively
to obtain
\begin{equation}
\mathrm{I}_\psi\leq\frac{t D_1}{2}\eps^2\delta^{-1}.
\label{U1}\end{equation}
Let $t\in[0,T]$, let $(t_n)\downarrow t$ and let $\nu^{i,+}_{x,t}$, be a weak-limit (in the sense of \refe{cvYoungMeasure}) of $\nu^{i,+}_{x,t_n}$. Then $\nu^{i,+}_{x,t}$ satisfies 
$$
\E\int_\R|\xi|^p d\nu^{i,+}_{x,t}(\xi) dx\leq C_p,
$$
and we have a similar bound for $\nu^{i,-}$. In particular, by \refe{e10}, $\chi_{f_i^\pm(t)}$ is integrable on $\T^N\times\R$ and 
\begin{multline}
\E\int_{\T^N}\int_{\R} f_1^\pm(x,t,\xi)\bar f_2^\pm(x,t,\xi) dx d\xi\nonumber\\
=\E\int_{(\T^N)^2}\int_{\R^2} \rho_\eps(x-y)\psi_\delta(\xi-\zeta)  f_1^\pm(x,t,\xi)\bar f_2^\pm(x,t,\xi)d\xi d\zeta dx dy +\eta_t(\eps,\delta),\label{U2}
\end{multline}
where $\lim_{\eps,\delta\to 0}\eta_t(\eps,\delta)=0$. 
To conclude, we need a bound on the term $\mathrm{I}_\rho$. Since $a$ has at most polynomial growth, there exists $C\geq 0$, $p> 1$, such that
\begin{equation*}
|a(\xi)-a(\zeta)|\leq \Gamma(\xi,\zeta)|\xi-\zeta|,\quad \Gamma(\xi,\zeta)=C(1+|\xi|^{p-1}+|\zeta|^{p-1}).
\end{equation*}
Supposing additionally that $\psi_\delta(\xi)=\delta^{-1}\psi_1(\delta^{-1}\xi)$ where $\psi_1$ is supported in $(-1,1)$, this gives
\begin{equation*}
|\mathrm{I}_\rho|\leq \E\int_0^t\int_{(\T^N)^2}\int_{\R^2}f_1\bar f_2\Gamma(\xi,\zeta)|\xi-\zeta|\psi_\delta(\xi-\zeta)|\nabla_x\rho_\eps(x-y)| d\xi d\zeta dx dy d\sigma.
\end{equation*}
By integration by parts with respect to $(\xi,\zeta)$, we deduce
\begin{equation*}
|\mathrm{I}_\rho|\leq \E\int_0^t\int_{(\T^N)^2}\int_{\R^2}\Upsilon(\xi,\zeta)d\nu^1_{x,\sigma}\otimes\nu^2_{y,\sigma}(\xi,\zeta)|\nabla_x\rho_\eps(x-y)|  dx dy d\sigma,
\end{equation*}
where
\begin{equation*}
\Upsilon(\xi,\zeta)=\int_\zeta^{+\infty}\int_{-\infty}^\xi\Gamma(\xi',\zeta')|\xi'-\zeta'|\psi_\delta(\xi'-\zeta') d\xi'd\zeta'.
\end{equation*}
It is shown below that $\Upsilon$ admits the bound
\begin{equation}
\Upsilon(\xi,\zeta)\leq C(1+|\xi|^{p}+|\zeta|^{p})\delta.
\label{ffp}\end{equation}
Since $\nu^1$ and $\nu^2$ vanish at infinity, we then obtain, for a given constant $C_p$,
\begin{equation*}
|\mathrm{I}_\rho|\leq t C_p\delta \left(\int_{\T^N}|\nabla_x\rho_\eps(x)|dx\right).
\end{equation*}
It follows that, for possibly a different $C_p$,
\begin{equation}
|\mathrm{I}_\rho|\leq t C_p  \delta\eps^{-1}.
\label{U3}\end{equation}
We then gather \refe{U1}, \refe{U3} and \refe{CR0} to deduce for $t\in [0,T]$
\begin{equation}
\E\int_{\T^N}\int_{\R} f_1^\pm(t)\bar f_2^\pm(t) dx d\xi\leq \int_{\T^N}\int_{\R} f_{1,0}\bar f_{2,0}dx d\xi+r(\eps,\delta),
\label{U4}\end{equation}
where the remainder $r(\eps,\delta)$ is
$
\ds r(\eps,\delta)=T C_p\delta\eps^{-1}+\frac{T D_1}{2}\eps^2\delta^{-1}+\eta_t(\eps,\delta)+\eta_0(\eps,\delta).
$
Taking $\delta=\eps^{4/3}$ and letting $\eps\to 0$ gives 
\begin{equation}
\E\int_{\T^N}\int_{\R} f_1^\pm(t)\bar f_2^\pm(t) dx d\xi\leq\int_{\T^N}\int_{\R} f_{1,0}\bar f_{2,0} dx d\xi. 
\label{L1compadd0}\end{equation}
Assume that $f$ is a generalized solution to \refe{stoSCL} with initial datum $\mathbf{1}_{u_0>\xi}$. Since
 $f_0$ is the (translated) Heaviside function $\mathbf{1}_{u_0>\xi}$, we have the identity $f_0\bar 
 f_0=0$. Taking $f_1=f_2=f$ in \refe{L1compadd0}, we deduce $f^+(1-f^+)=0$ a.e., {\it i.e.} 
$f^+\in\{0,1\}$ a.e. The fact that $-\partial_\xi f^+$ is a Young measure then gives the conclusion: 
indeed, by Fubini Theorem, for any $t\in [0,T]$, there is a set $E_t$ of full measure in $\T^N\times
\Omega$ such that, for  $(x,\omega)\in E_t$, $f^+(x,t,\xi,\omega)\in\{0,1\}$ for a.e. $\xi\in\R$. Recall 
that $-\partial_\xi f^+(x,t,\cdot,\omega)$ is a probability measure on $\R$ so that, necessarily, there 
exists $u^+(x,t,\omega)\in\R$ such that $f^+(t,x,\xi,\omega)=\mathbf{1}_{u^+(x,t,\omega)>\xi}$ 
for almost every $(x,\xi,\omega)$. 
In particular, $u^+=\int_\R (f^+-\mathbf{1}_{\xi>0}) d\xi$ for almost every $(x,\omega)$. We have a similar result for $f^-$.

The discussion 
after Definition \ref{d4} tells us that $f^+$ being solution in the sense of Definition \ref{d4} implies that $u^+$ is a solution in the sense of Definition \ref{defkineticsol}. Since $f=f^+$ a.e., this shows the reduction of generalized solutions to solutions. If now $u_1$ and $u_2$ are two solutions to \refe{stoSCL}, we deduce from \refe{L1compadd0} with $f_i=\mathbf{1}_{u_i>\xi}$ and from the identity
\begin{equation*}
\int_\R \mathbf{1}_{u_1>\xi}\overline{\mathbf{1}_{u_2>\xi}}d\xi=(u_1-u_2)^+
\end{equation*}
the contraction property
\begin{equation}
\E\|(u_1(t)-u_2(t))^+\|_{L^1(\T^N)}\leq\E\|(u_{1,0}-u_{2,0})^+\|_{L^1(\T^N)}.
\label{L1compadd}\end{equation}
This implies the $L^1$-contraction property, comparison and uniqueness of solutions.
\medskip

In the multiplicative case ($\Phi$ depending on $u$), the reasoning is similar, except that there is an additional term in the bound on $\mathrm{I_\psi}$. More precisely, by Hypothesis~\refe{D1} we obtain in place of \refe{U1} the estimate 
\begin{equation*}
\ds\mathrm{I_\psi}\leq \frac{T D_1}{2}\eps^2\delta^{-1}
+\frac{D_1}{2}\mathrm{I}^h_\psi,
\end{equation*} 
where
\begin{equation*}
\mathrm{I}_\psi^h=\E\int_0^t\int_{(\T^N)^2}\rho_\eps\int_{\R^2} \psi_\delta(\xi-\zeta)|\xi-\zeta|h(|\xi-\zeta|)d\nu^1_{x,\sigma}\otimes\nu^2_{y,\sigma}(\xi,\zeta)dx dy d\sigma.
\end{equation*}
Choosing $\psi_\delta(\xi)=\delta^{-1}\psi_1(\delta^{-1}\xi)$ with $\psi_1$ compactly supported gives
\begin{equation}
\mathrm{I_\psi}\leq \frac{T D_1}{2}\eps^2\delta^{-1}+\frac{TD_1C_\psi h(\delta)}{2},\quad C_\psi:=\sup_{\xi\in\R}\|\xi\psi_1(\xi)\|.
\label{U1bis}\end{equation}
We deduce \refe{U4} with a remainder term $\ds r'(\eps,\delta):=r(\eps,\delta)+\frac{TD_1C_\psi h(\delta)}{2}$ and conclude the proof as in the additive case.
\medskip

There remains to prove \refe{ffp}: setting $\xi''=\xi'-\zeta'$, we have
\begin{align*}
\Upsilon(\xi,\zeta)=&\int_\zeta^{+\infty}\int_{|\xi''|<\delta,\xi''<\xi-\zeta'}\Gamma(\xi''+\zeta',\zeta')|\xi''|\psi_\delta(\xi'') d\xi''d\zeta'\\
\leq&C\int_{\zeta}^{\xi+\delta}\max_{|\xi''|<\delta,\xi''<\xi-\zeta'}\Gamma(\xi''+\zeta',\zeta') d\zeta'\; \delta\\
\leq &C \int_{\zeta}^{\xi+\delta}(1+|\xi|^{p-1}+|\zeta'|^{p-1}) d\zeta'\; \delta,
\end{align*}
which gives \refe{ffp}. \qed
\bigskip

{\bf Proof of Corollary~\ref{cor:timecontinuity}:} In the proof of Theorem~\ref{th:Uadd}, we have shown that there exists $u^+$ such that for every $t\in[0,T)$, for almost all $\omega,x,\xi$, $f^+(\omega,x,t,\xi)=\mathbf{1}_{u^+(\omega,x,t)>\xi}$. We will show that $u^+$ has almost surely con\-ti\-nu\-ous trajectories. Since $u=u^+$ a.e. with respect to $(\omega,t,x)$, this will give the result. Let us first prove that, a.s., $u^+$ has left and right limits at every $t\in(0,T)$. With a similar proof, we will obtain that $u^+$ also has a right limit at $t=0$. By Remark~\ref{rk:integrabilityfpm}, and by considering an increasing sequence of exponent $p$, we can fix $\omega$ in a set of full measure such that 
\begin{equation}\label{eq:integrabilityuomega}
\sup_{t\in[0,T]}\|u^+(t)\|_{L^p(\T^N)}\leq C_p(\omega)
\end{equation}
for all $1\leq p<+\infty$. Let $t\in(0,T)$ and let $(t_n)\downarrow t$. By Proposition~\ref{prop:LRlimits} applied to the solution $f^+$, the weak-star limit of $f^+(t_n)$ exists in $L^\infty(\Omega\times\T^N\times\R)$. This limit is also the limit of 
$$
\frac{1}{\eps}\int_t^{t+\eps}f^+(s)ds=\frac{1}{\eps}\int_t^{t+\eps}f(s)ds
$$
as $\eps\to 0$. It is therefore $f^+(t)$. Since $f^+(t)=\mathbf{1}_{u^+(t)>\xi}$ is at equilibrium, Lemma~\ref{lem:weakstrongeq} and \refe{eq:integrabilityuomega} give $u^+(t_n)\to u^+(t)$ in $L^p(\T^N)$. Similarly, we use the fact that $f^-$ is at equilibrium to prove the existence of a left limit. Let us now first show the continuity at $t=0$: this is a consequence of Remark~\ref{RKequilibrium}, we have $f^{+,0}=\mathbf{1}_{u_0>\xi}$. In particular,
$$
u^+(x,0)=\int_\R (f^{0,+}(x,\xi)-\mathbf{1}_{0>\xi}) d\xi=\int_\R (\mathbf{1}_{u_0(x)>\xi}-\mathbf{1}_{0>\xi}) d\xi=u_0(x).
$$
To prove similar results at time $t_*\in(0,T)$, we consider $t_*$ as the origin of time: indeed it follows from \refe{eq:kineticfpre} and Proposition~\ref{prop:LRlimits} that 
\begin{multline*}
\int_{t_*}^T\<f^+(t),\partial_t \varphi(t)\>dt+\<f^{-}(t_*),\varphi(t_*)\>
+\int_{t_*}^T \<f^+(t),a(\xi)\cdot\nabla\varphi(t)\>dt\\
=-\sum_{k\geq 1}\int_{t_*}^T\int_{\T^N}\int_\R g_k(x,\xi)\varphi(x,t,\xi)d\nu^+_{x,t}(\xi)dxd\beta_k(t)\\
-\frac{1}{2}\int_{t_*}^T\int_{\T^N}\int_\R \partial_\xi\varphi(x,t,\xi)\GG^2(x,\xi)d\nu^+_{(x,t)}(\xi) dx dt+m(\mathbf{1}_{[t_*,T]}\partial_\xi\varphi).
\end{multline*}
In other words, $t\mapsto f^+(t_*+t)$ is a generalized solution to \refe{stoSCL} on $[0,T-t_*]$ with initial datum $f^-(t_*)=\mathbf{1}_{u^-(t_*)>\xi}$. We obtain $u^+(t_*)=u^-(t_*)$ and the 
result follows. \qed

\subsection{Entropy solutions}\label{s3.3}

For deterministic first-order scalar conservation laws, the notion of entropy solution was introduced by Kruzhkov~\cite{Kruzhkov70} prior to the notion of kinetic solution~\cite{LionsPerthameTadmor94}. For the first-order scalar conservation law with stochastic forcing, a corresponding notion of weak entropy solution has been introduced by Feng and Nualart~\cite{FengNualart08}:

\begin{definition}[Weak entropy solution] A measurable function $u\colon\T^N\times [0,T]\times\Omega\to\R$ is said to be a weak entropy solution to~\refe{stoSCL}  if $(u(t))$ is an adapted $L^2(\T^N)$-valued process, for all $p\geq 1$, there exists $C_p\geq 0$ such that 
\begin{equation*}
\E\left(\esssup_{t\in[0,T]}\|u(t)\|_{L^p(\T^N)}^p\right)\leq C_p,
\end{equation*} 
and for all convex $\eta\in C^2(\R)$, for all non-negative $\theta\in C^1(\T^N)$, for all $0\leq s\leq t\leq T$,
\begin{multline}
\<\eta(u(t)),\theta\>-\<\eta(u(s)),\theta\>
\leq
\int_s^t \<q(u(r)),\nabla\theta\>dr\\
+\sum_{k\geq 1}\int_s^t\<g_k(\cdot,u(r))\eta'(u(r)),\theta\>d\beta_k(r)
+\frac{1}{2}\int_s^t \<\GG^2(\cdot,u(r))\eta''(u(r)),\theta\> dr,
\label{eq:weakentropy}\end{multline}
a.s., where $q(u)=\int_0^u a(\xi)\eta'(\xi)d\xi$.
\label{def:entropysol}\end{definition}

An entropy solution is a kinetic solution and vice versa. To prove this fact, let us introduce an auxiliary definition:

\begin{definition}[Time-weak weak entropy solution] Let $u_0\in L^\infty(\T^N)$. A measurable function $u\colon\T^N\times [0,T]\times\Omega\to\R$ is said to be a time-weak weak entropy solution to~\refe{stoSCL} with initial datum $u_0$ if $(u(t))$ is an adapted $L^2(\T^N)$-valued process, for all $p\geq 1$, there exists $C_p\geq 0$ such that 
\begin{equation*}
\E\left(\esssup_{t\in[0,T]}\|u(t)\|_{L^p(\T^N)}^p\right)\leq C_p,
\end{equation*} 
and for all convex $\eta\in C^2(\R)$, for all non-negative $\rho\in C^1_c(\T^N\times[0,T))$, 
\begin{multline}
\int_0^T \<\eta(u),\partial_t\rho\>dr+\<\eta(u_0),\rho(0)\>-
\int_0^T \<q(u),\nabla\rho\>dr\\\geq
-\sum_{k\geq 1}\int_0^T\<g_k(\cdot,u(r))\eta'(u(r)),\rho\>d\beta_k(r)
-\frac{1}{2}\int_0^T \<\GG^2(\cdot,u(r))\eta''(u(r)),\rho\> dr,
\label{eq:timeweakentropy}\end{multline}
a.s., where $q(u)=\int_0^u a(\xi)\eta'(\xi)d\xi$.
\label{def:timeweakentropysol}\end{definition}

\begin{proposition}[Entropy and kinetic solutions] Let $u_0\in L^\infty(\T^N)$. For a measurable function $u\colon\T^N\times [0,T]\times\Omega\to\R$, it is equivalent to be a kinetic solution to \refe{stoSCL}, {\it i.e.} a solution in the sense of Definition~\ref{defkineticsol}, and a time-weak weak solution.
\label{prop:EQUentropykinetic}\end{proposition}

The proof of the proposition is classical. Choosing test functions $\varphi(x,t,\xi)=\rho(x,t)\eta'(\xi)$ in \refe{eq:kineticupre} and using the inequality $m\eta''\geq 0$ gives \refe{eq:timeweakentropy}. Conversely, starting from \refe{eq:timeweakentropy}, one defines the measure $m$ (actually $\partial^2_\xi m$) by 
\begin{multline*}
m(\rho\otimes\eta'')=\int_0^T \<\eta(u),\partial_t\rho\>dr+\<\eta(u_0),\rho(0)\>+
\int_0^T \<q(u),\nabla\rho\>dr\\
+\sum_{k\geq 1}\int_0^T\<g_k(\cdot,u(r))\eta'(u(r)),\rho\>d\beta_k(r)
+\frac{1}{2}\int_0^T \<\GG^2(\cdot,u(r))\eta''(u(r)),\rho\> dr,
\end{multline*}
and then derives \refe{eq:kineticupre}. See \cite{PerthameBook} for precise references.
\medskip

It is clear also that a weak entropy solution, satisfying $u(0)=u_0$, is a time-weak entropy solution, 
while, for the converse assertion, time-continuity of the solution is required. We have seen that a 
kinetic solution is continuous in time, it follows that it is indeed a weak entropy solution. 
Then, taking $\eta(u) = |u|^p$ and $\theta=1$ in \eqref{eq:weakentropy}, it is classical to prove that \eqref{eq:integrabilityu} is satisfied.

\subsection{Spatial regularity}

To conclude this paragraph and our applications of Proposition~\ref{prop:CR}, we give a result on the spatial regularity of the solution. To that purpose, we introduce two semi-norms that measure the $W^{\sigma,1}$-regularity of a function $u\in L^1(\T^N)$ ($\sigma\in(0,1)$): we set 
\begin{equation*}
p^\sigma(u):=\int\limits_{\T^N}\int\limits_{\T^N}\frac{|u(x)-u(y)|}{|x-y|^{N+\sigma}}dx dy,
\end{equation*}
and
\begin{equation*}
p^\sigma_\rho(u)=\sup_{0<\eps<2D_N}\frac{1}{\eps^\sigma}\int\limits_{\T^N}\int\limits_{\T^N}|u(x)-u(y)|\rho_\eps(x-y)dx dy,
\end{equation*}
where $(\rho_\eps)$ is a fixed regularizing kernel: $\rho_\eps(x)=\eps^{-N}\rho(\eps^{-1}|x|)$ where $\rho$ is supported in the ball $B(0,1)$ of $\R^N$ and where $D_N=\sqrt{N}$ is the diameter of $[0,1]^N$. We define $W^{\sigma,1}(\T^N)$ as the subspace of $u\in L^1(\T^N)$ with finite norm
\begin{equation*}
\|u\|_{W^{\sigma,1}(\T^N)}=\|u\|_{L^1(\T^N)}+p^\sigma(u).
\end{equation*}

\begin{lemma}[Comparison of the $W^{\sigma,1}$ semi-norms] Let $\sigma\in(0,1)$. Then there exists $C$ depending on $\sigma$, $\rho$, $N$ such that, for all $0<s<\sigma$, for all $u\in L^1(\T^N)$,
\begin{equation*}
p^\sigma_\rho(u)\leq C p^\sigma(u),\quad p^s(u)\leq \frac{C}{\sigma-s}\; p^\sigma_\rho(u).
\end{equation*}
\label{lem:comparisonWsigma}\end{lemma}

{\bf Proof:} we have 
\begin{equation*}
\frac{1}{\eps^\sigma}\rho_\eps(x-y)\leq \frac{\|\rho\|_{L^\infty}}{\eps^{N+\sigma}}\mathbf{1}_{|x-y|<\eps}\leq\frac{\|\rho\|_{L^\infty}}{|x-y|^{N+\sigma}},
\end{equation*}
hence $p^\sigma_\rho(u)\leq C p^\sigma(u)$. Conversely, we multiply the inequality
\begin{equation*}
\frac{1}{\eps^\sigma}\int\limits_{\T^N}\int\limits_{\T^N}|u(x)-u(y)|\frac{1}{\eps^N}\rho\left(\frac{|x-y|}{\eps}\right)dx dy\leq p^\sigma_\rho(u)
\end{equation*}
by $\eps^{-1+(\sigma-s)}$ and sum over $\eps\in(0,2D_N)$. We obtain 
$$
\alpha p^s(u)\leq \frac{C}{\sigma-s}\; p^\sigma_\rho(u),
$$
where
$$
\alpha:=\int_{1/2}^1 \tau^{s+N-1}\rho(\tau)d\tau>0,
$$
which gives the second inequality. \qed

\begin{remark}[A third semi-norm] The proof of Lemma~\ref{lem:comparisonWsigma} involves a third $W^{s,1}$-semi-norm, that we will actually prefer to the first two introduced above. Indeed it shows that, for some given constant $C,\alpha>0$, we have, for all $0<s<\sigma$, for all $u\in L^1(\T^N)$, 
$$
\alpha p^s(u)\leq \widetilde{p^s_\rho}(u)\leq \frac{C}{\sigma-s}\; p^\sigma_\rho(u),
$$
where
$$ 
\widetilde{p^s_\rho}(u):=\int_0^{2D_N} \frac{1}{\eps^{s+1}}\int\limits_{\T^N}\int\limits_{\T^N}|u(x)-u(y)|\rho_\eps(x-y)dx dy d\eps.
$$
\label{rk:defprho}\end{remark}

\begin{theorem}[$W^{\sigma,1}$-regularity] Let $u_{0}\in L^\infty(\T^N)$, let $u\colon\T^N\times(0,+\infty)\times\Omega\to\R$ be the solution to~\refe{stoSCL} with initial datum $u_0$. Assume that $h$ satisfies 
\begin{equation}
h(\delta)\leq C \delta^\alpha,\quad \delta<1,\quad 0<\alpha.
\label{decreasingh}\end{equation}
Set $\sigma=\min\left(\frac{2\alpha}{1+\alpha},\frac{1}{2}\right)$. Then, there exists a constant $C$ such that, for all $t\geq 0$, we have
\begin{equation}
\E \widetilde p^{\sigma}_\rho(u(t))\leq C( p_\rho^{\sigma}(u_0)+t).
\label{eq:Wu}\end{equation}
In particular, for all $0<s<\sigma$, there exists a constant $C_s>0$ such that for $t\ge 0$, 
$$
\E \|u(t)\|_{W^{s,1}(\T^N)}\leq C_s(\|u_0\|_{W^{\sigma,1}(\T^N)}+t).
$$
\label{th:spaceregularity}\end{theorem}

{\bf Proof:} the last assertion is proved as follows: by Lemma~\ref{lem:comparisonWsigma}, \refe{eq:Wu} implies $\E p^s(u(t))\leq C_s(p^{\sigma}(u_0)+t)$. Poincar\'e Inequality gives
\begin{equation*}
\left\|u(t)-\int_{\T^N}u(t)dx\right\|_{L^1(\T^N)}\leq C_s p^s(u(t)).
\end{equation*}
Since $\E\int_{\T^N}u(t)dx=\E\int_{\T^N}u_0 dx$, we obtain a bound on the $L^1$-norm of $u$:
\begin{equation*}
\left\|u(t)\right\|_{L^1(\T^N)}\leq C_s (p^{s}(u_0)+t+\|u_0\|_{L^1(\T^N)}),
\end{equation*}
hence $\E \|u(t)\|_{W^{s,1}(\T^N)}\leq C_s(\|u_0\|_{W^{\sigma,1}(\T^N)}+t)$. To prove \refe{eq:Wu}, we apply Prop.~\ref{prop:CR} with $f_1=f_2=\mathbf{1}_{u>\xi}$, $\rho=\rho^\eps$, $\psi=\psi_\delta$. Since $\partial_\xi\mathbf{1}_{u>\xi}=-\delta_{u=\xi}$ is a Radon measure with mass $1$, we have
\begin{multline*}
\E\int_{(\T^N)^2}\rho_\eps(x-y)(u(x,t)-u(y,t))^+ dx dy\\
\leq\E\int_{(\T^N)^2}\int_{\R^2} \rho_\eps(x-y)\psi_\delta(\xi-\zeta)\mathbf{1}_{u(x,t)>\xi}(1-\mathbf{1}_{u(y,t)>\zeta}) dx dy d\xi d\zeta+\delta
\end{multline*}
and
\begin{multline*}
\E\int_{(\T^N)^2}\int_{\R^2} \rho_\eps(x-y)\psi_\delta(\xi-\zeta)\mathbf{1}_{u_0(x)>\xi}(1-\mathbf{1}_{u_0(y)>\zeta}) dx dy d\xi d\zeta\\
\leq\E\int_{(\T^N)^2}\rho_\eps(x-y)(u_0(x)-u_0(y))^+ dx dy+\delta.
\end{multline*}
We deduce that
\begin{multline*}
\E\int_{(\T^N)^2}\rho_\eps(x-y)(u(x,t)-u(y,t))^+ dx dy\\
\leq\E\int_{(\T^N)^2}\rho_\eps(x-y)(u_0(x)-u_0(y))^+ dx dy+\mathrm{I}_\rho+\mathrm{I}_\psi+2\delta.
\end{multline*}
As in \refe{U1bis}-\refe{U3}, we have  
\begin{equation*}
\mathrm{I}_\psi\leq t C(\eps^2\delta^{-1}+h(\delta)),\quad \mathrm{I}_\rho\leq t C  \delta\eps^{-1},
\end{equation*}
hence
\begin{multline*}
\E\int_{(\T^N)^2}\rho_\eps(x-y)(u(x,t)-u(y,t))^+ dx dy\\
\leq\E\int_{(\T^N)^2}\rho_\eps(x-y)(u_0(x)-u_0(y))^+ dx dy+ t C(\eps^2\delta^{-1}+h(\delta)+ \delta\eps^{-1})+2\delta.
\end{multline*}
By optimization in $\delta$, using \refe{decreasingh} and Remark~\ref{rk:defprho}, we obtain \refe{eq:Wu}. \qed

\section{Existence}\label{sec:exists}

\subsection{The parabolic approximation, kinetic formulation}\label{sec:viscousappx}

Let $u_{0}\in L^\infty(\T^N)$. To prove the existence of a solution to \refe{stoSCL} with initial datum $u_0$, we show the convergence of the parabolic approximation
\begin{equation}\left\{\begin{array}{r l l}
du^\eta+\div(A^\eta(u^\eta))dt-\eta\Delta u^\eta dt&=&\Phi_\eta(u^\eta)dW(t),\quad t>0,x\in\T^N,\\
u^\eta(x,0)&=&u_0^\eta(x),\quad x\in\T^N.
\end{array}\right.\label{stoSCL2}\end{equation}
Where $u_0^\eta$ is a smooth approximation of $u_0$,  $\Phi_\eta$ is a suitable Lipschitz approximation of $\Phi$ satisfying \eqref{D0}, \eqref{D1} uniformly. We define $g^\eta_k$ and $\GG^\eta$ as in the case $\eta=0$. 
It is possible to choose $g_k^\eta$ with compact support and smooth with respect to 
$(x,u)$. Moreover, we may assume that $g_k=0$ for $k\ge \frac1\eta$. Finally, we choose $A^\eta$ 
which is a smooth approximation of $A$ and has the same growth as $A$. We set $a^\eta=(A^\eta)'$.
 
It is shown in \cite{GyongyRovira00} that equation \eqref{stoSCL2} has a unique $L^\rho(\T^N)$ valued con\-ti\-nu\-ous solution provided $\rho$ is large enough and $u_0\in L^\rho(\T^N)$, hence in particular for 
 $u_{0}\in L^\infty(\T^N)$. 
Moreover, it is also shown in \cite{GyongyRovira00} that using  It\^o Formula one can prove 
that $u^\eta$ satisfies the energy inequality
\begin{multline}
\E\|u^\eta(t)\|^2_{L^2(\T^N)}+2\eta\E\int_0^t \|\nabla u^\eta\|_{L^2(\T^N)}^2 ds\\
\leq \E\|u_0^\eta\|^2_{L^2(\T^N)}+\E\int_0^t \|\GG_\eta(u^\eta)\|_{L^2(\T^N)}^2 ds.
\label{energyineq}\end{multline}
By \eqref{D0} and Gronwall Lemma, we easily derive 
\begin{equation}
\E\|u^\eta(t)\|^2_{L^2(\T^N)}+\eta\E\int_0^t \|\nabla u^\eta\|_{L^2(\T^N)}^2 ds\leq C(T)( \E\|u_0^\eta\|^2_{L^2(\T^N)}+1).
\label{e32}\end{equation}
Also, for $p\ge 2$, by It\^o Formula applied to $|u|^p$ and a martingale inequality
\begin{equation}
\E\left(\sup_{t\in[0,T]}\|u^\eta(t)\|^p_{L^p(\T^N)}\right)+\eta\E\int_0^T\int_{\T^N} |u^\eta(t,x)|^{p-2}|\nabla u^\eta(t)|^2 dxdt \leq C(p,u_0^\eta,T).
\label{estimLp}\end{equation}


\begin{proposition}[Kinetic formulation] Let $u_{0}^\eta\in C^3(\T^N)$ and let $u^\eta$ be the solution to \refe{stoSCL2}. Then $f^\eta:=\mathbf{1}_{u^\eta>\xi}$ satisfies: for all $\varphi\in C^1_c(\T^N\times[0,T)\times\R)$, 
\begin{multline}
\int_0^T\<f^\eta(t),\partial_t \varphi(t)\>dt+\<f_0,\varphi(0)\>
+\int_0^T \<f^\eta(t),a^\eta(\xi)\cdot\nabla\varphi(t)-\eta\Delta\varphi(t)\>dt\\
=-\sum_{k\geq 1}\int_0^T\int_{\T^N}\int_\R g_k^\eta(x,\xi)\varphi(x,t,\xi)d\nu^\eta_{x,t}(\xi)dxd\beta_k(t)\\
-\frac{1}{2}\int_0^T\int_{\T^N}\int_\R \partial_\xi\varphi(x,t,\xi)\GG^2_\eta(x,\xi)d\nu^\eta_{(x,t)}(\xi) dx dt+m^\eta(\partial_\xi\varphi),
\label{eq:kineticueta}\end{multline}
a.s., where $f_0(\xi)=\mathbf{1}_{u_0>\xi}$ and, for $\phi\in C_b(\T^N\times[0,T]\times\R)$, 
\begin{equation*}
\nu^\eta_{(x,t)}=\delta_{u^\eta(x,t)},\quad m^\eta(\phi)=\int_{\T^N\times[0,T]\times\R} \phi(x,t,u^\eta(x,t))\eta |\nabla u^\eta|^2 dx dt.
\end{equation*}
\end{proposition}

Note that the measure $m^\eta$ is explicitly known here: $m^\eta=\eta |\nabla u^\eta|^2\delta_{u^\eta=\xi}$.
\medskip

{\bf Proof:} By \cite{hofmanova-nodea}, we know that $u^\eta$ is almost 
surely continuous with values in $C^2(\T^N)$. Thus we may use It\^o Formula and obtain, for $\theta\in C^2(\R)$ with polynomial growth at $\pm\infty$,
\begin{align*}
d(\mathbf{1}_{u^\eta>\xi},\theta')& :=d\int_\R \mathbf{1}_{u^\eta>\xi}\theta'(\xi) d\xi=d\theta(u^\eta)\\
	&=\theta'(u^\eta)(-a^\eta(u^\eta)\cdot\nabla u^\eta dt+\eta\Delta u^\eta dt+\Phi_\eta(u^\eta) dW)+\frac{1}{2}\theta''(u^\eta)\GG_\eta^2 dt.
\end{align*}
We rewrite the first term as
\begin{equation*}
-\theta'(u^\eta)a^\eta(u^\eta)\cdot\nabla u^\eta =-\div\left\{\int_0^{u^\eta} a^\eta(\xi)\theta'(\xi)d\xi\right\}
=-\div(a^\eta\mathbf{1}_{u^\eta>\xi},\theta'),
\end{equation*}
the second term as
\begin{equation*}
\begin{array}{ll}
\ds \theta'(u^\eta)\eta\Delta u^\eta &\ds =\eta\Delta\theta(u^\eta)dt-\eta|\nabla u^\eta|^2\theta''(u^\eta)\\
\\
&\ds =
\eta\Delta(\mathbf{1}_{u^\eta>\xi},\theta')
+(\partial_\xi(\eta|\nabla u^\eta|^2\delta_{u^\eta=\xi}),\theta')
\end{array}
\end{equation*}
to obtain the kinetic formulation
\begin{multline}
d(\mathbf{1}_{u^\eta>\xi},\theta')=-\div[(a^\eta\mathbf{1}_{u^\eta>\xi},\theta')]dt+\eta\Delta(\mathbf{1}_{u^\eta>\xi},\theta')dt\\
+(\partial_\xi(\eta|\nabla u^\eta|^2\delta_{u^\eta=\xi}-\frac{1}{2}\GG_\eta^2\delta_{u^\eta=\xi}),\theta')dt+\sum_{k\geq 1}(\delta_{u^\eta=\xi},\theta' g_{k,\eta})d\beta_k.
\label{kineticeta0}\end{multline}
Taking $\theta(\xi)=\int_{-\infty}^\xi\beta$, we then obtain \refe{kineticeta0} with the test function 
$\beta$ in place of $\theta'$. Since the test functions $\varphi(x,\xi)=\alpha(x)\beta(\xi)$ form a 
dense subset of $C^\infty_c(\T^N\times\R)$, \refe{eq:kineticueta} follows. \qed
\medskip

Equation~\refe{eq:kineticueta} is close to the kinetic equation \refe{eq:kineticupre} satisfied by the solution to \refe{stoSCL}. For $\eta\to 0$, we lose the precise structure of $m^\eta=\eta|\nabla u^\eta|^2\delta_{u^\eta=\xi}$ and obtain a solution $u$ to \refe{stoSCL}. More precisely, we will prove the
\begin{theorem}[Convergence of the parabolic approximation] Let $u_{0}\in L^\infty(\T^N)$. There exists a unique solution $u$ to \refe{stoSCL} with initial datum $u_0$ which is the strong limit of $(u^\eta)$ as $\eta\to0$: for every $T>0$, for every $1\leq p<+\infty$,
\begin{equation}
\ds\lim_{\eta\to0}\E\|u^\eta-u\|_{L^p(\T^N\times(0,T))}=0.
\label{eqcvexist}\end{equation}
\label{th:cvexists}\end{theorem}

The proof of Theorem~\ref{th:cvexists} is quite a straightforward consequence of both the result of reduction of generalized solution to solution - Theorem~\ref{th:Uadd} - and the a priori estimates derived in the following section.

\subsubsection{A priori estimates}\label{sec:apriorieta}

We denote indifferently by $C_p$ various constants that may depend on $p\in[1,+\infty)$, on $u_0$, on the noise and on the terminal time $T$, but not on $\eta\in(0,1)$.
\medskip

{\bf 1. Estimate of $m^\eta$:}
we analyze the kinetic measure $m^\eta=\eta|\nabla u^\eta|^2\delta_{u^\eta=\xi}$.  
By 
\eqref{e32},  we have a uniform bound $\E m^\eta(\T^N\times[0,T]\times\R)\leq C$. Furthermore, the second term in the left hand-side of \eqref{estimLp} is $\E\int_{\T^N\times[0,T]\times\R} |\xi|^{p-2}dm^\eta(x,t,\xi)$, so we have
\begin{equation}
\E\int_{\T^N\times[0,T]\times\R} |\xi|^{p}dm^\eta(x,t,\xi)\leq C_p.
\label{estimmeta0}\end{equation}
We also have the improved estimate, for $p\geq 0$,
\begin{equation}
\E\left|\int_{\T^N\times[0,T]\times\R} |\xi|^{2p}dm^\eta(x,t,\xi)\right|^2\leq C_p.
\label{estimmeta}\end{equation}
To prove \refe{estimmeta}, we apply It\^o Formula to $\psi(u^\eta)$, $\psi(\xi):=|\xi|^{2p+2}$: 
\begin{equation*}
d\psi(u^\eta)+\div(\mathcal{F})dt+\eta\psi''(u^\eta)|\nabla u^\eta|^2 dt
=\psi'(u^\eta)\Phi_\eta(u^\eta)dW+\frac12\psi''(u^\eta)\GG_\eta^2 dt,
\end{equation*}
where $\mathcal{F}:=\int_0^{u^\eta}a^\eta(\xi)\psi'(\xi)d\xi-\eta\nabla \psi(u^\eta)$.
It follows
\begin{multline*}
\int_0^T\int_{\T^N}\eta\psi''(u^\eta)|\nabla u^\eta|^2 dx dt\\
\leq \int_{\T^N}\psi(u_0) dx
+\sum_{k\geq 1}\int_0^T\int_{\T^N}\psi'(u^\eta)g_{k,\eta}(x,u^\eta) dx d\beta_k(t)\\
+\frac12\int_0^T\int_{\T^N}\GG_\eta^2(x,u^\eta)\psi''(u^\eta)dx dt.
\label{itoeta01}\end{multline*}
Taking the square, then expectation, we deduce by It\^o isometry
\begin{multline*}
\E\left|\int_0^T\int_{\T^N}\eta\psi''(u^\eta)|\nabla u^\eta|^2 dx dt\right|^2 \leq 3\E\left|\int_{\T^N}\psi(u_0) dx\right|^2\\
+3\E\int_0^T\sum_{k\geq 1}\left|\int_{\T^N} g_k(x,u^\eta)\psi'(u^\eta)dx \right|^2 dt
+\frac32\E\left|\int_0^T\int_{\T^N}\GG^2(x,u^\eta)\psi''(u^\eta)dx dt\right|^2.
\end{multline*}
By \refe{D0}, \refe{estimLp} and Cauchy-Schwarz inequality, we obtain \refe{estimmeta}.

\medskip

{\bf 2. Estimate on $\nu^\eta$:} By the bound~\refe{estimLp} on $u^\eta$ in $L^p$, we have,
\begin{equation}
\E\esssup_{t\in[0,T]}\int_{\T^N}\int_\R |\xi|^p d\nu^\eta_{x,t}(\xi)dx\leq C_p
\label{boundLpnueta}\end{equation}
and, in particular,
\begin{equation}
\E\int_0^T\int_{\T^N}\int_\R |\xi|^p d\nu^\eta_{x,t}(\xi)dx dt\leq C_p.
\label{boundLpnuetaT}\end{equation}

\subsubsection{Generalized solution}\label{sec:generalizedsoleta}
Consider a sequence $(\eta_n)\downarrow 0$. We use the a priori bounds derived in the preceding subsection to deduce, up to subsequences:
\begin{enumerate}
\item by \refe{boundLpnuetaT} and Theorem~\ref{th:youngmeasure} and 
Corollary~\ref{cor:kineticfunctions} respectively, the convergence $\nu^{\eta_n}\to\nu$ (in the sense of 
\refe{cvYoungMeasure}) and the convergence $f^{\eta_n}\rightharpoonup f$ in $L^\infty(\Omega\times
\T^N\times(0,T)\times\R)$-weak-*. Besides, the bound \refe{boundLpnueta} is stable: $\nu$ satisfies 
\refe{eq:integrabilityf}.
\item 

For $r\in\N^*$, let $K_r=\T^N\times[0,T]\times[-r,r]$ and let $\mathcal{M}_r$ denote the space of bounded Borel measures over $K_r$ (with norm given by the total variation of measures). It is the topological dual of $C(K_r)$, the set of 
continuous functions on $K_p$. Since $\mathcal{M}_r$ is separable, 
the space $L^2(\Omega;\mathcal{M}_r)$ is 
the topological dual space of $L^2(\Omega,C(K_r))$, {\it c.f.} Th\'eor\`eme~1.4.1 in \cite{Droniou01}. The estimate \refe{estimmeta} with $p=0$ gives a uniform bound on $(m^{\eta_n})$ in $L^2(\Omega,\mathcal{M}_r)$: there exists $m_r\in L^2(\Omega,\mathcal{M}_r)$ such that  
up to subsequence, $m^{\eta_n}\rightharpoonup m$ in $L^2(\Omega;\mathcal{M}_r)$-weak star. By a dia\-go\-nal process, we obtain $m_r=m_{r+1}$ in $L^2(\Omega;\mathcal{M}_r)$ and the convergence in all the spaces $L^2(\Omega;\mathcal{M}_r)$-weak star of a single subsequence still denoted $(m^{\eta_n})$. The condition at infinity \refe{estimmeta} then shows that $m$ defines an element of $L^2(\Omega;\mathcal{M})$, where $\mathcal{M}$ denotes the space of bounded Borel measures over $\T^N\times[0,T]\times\R$, that $m$ satisfies all the points {\it 1.} and {\it 2.} of Definition~\ref{def:kineticmeasure} and that
\begin{multline}\label{cvmntom}
\E\left(\alpha\int_{\T^N\times[0,T]\times\R}\phi(x,t,\xi)dm^{\eta_n}(x,s,\xi)\right)\\
\to\E\left(\alpha\int_{\T^N\times[0,T]\times\R}\phi(x,t,\xi)dm(x,s,\xi)\right),
\end{multline}
for every $\alpha\in L^2(\Omega)$, $\phi\in C_b(\T^N\times[0,T]\times\R)$. Let us check point {\it 3.} of Definition~\ref{def:kineticmeasure}: let  $\phi\in C_b(\T^N\times\R)$ and set 
$$
x^n(t):=\int_{\T^N\times[0,t]\times\R}\phi(x,\xi)dm^{\eta_n}(x,s,\xi),
$$
where $\alpha\in L^2(\Omega)$, $\gamma\in L^2([0,T])$. By Fubini's Theorem,
\begin{equation*}
\E\left(\alpha \int_0^T \gamma(t) x^n(t) dt\right)=\E\left(\alpha\int_{\T^N\times[0,t]\times\R}\phi(x,\xi)\Gamma(s) dm^{\eta_n}(x,s,\xi)\right),
\end{equation*} 
where $\Gamma(s)=\int_s^T \gamma(t)dt$. Since $\Gamma$ is continuous, \refe{cvmntom} gives
\begin{equation*}
\E\left(\alpha \int_0^T \gamma(t) x^n(t) dt\right)\to \E\left(\alpha \int_0^T \gamma(t) x(t) dt\right), 
\end{equation*}
where
\begin{equation*}
x(t)=\int_{\T^N\times[0,t]\times\R}\phi(x,\xi)dm(x,s,\xi).
\end{equation*}
Since tensor functions are dense in $L^2(\Omega\times[0,T])$, we obtain the weak convergence $x_n\to x$ in $L^2(\Omega\times[0,T])$. In particular, since the space of predictable process is weakly-closed, $x$ is predictable.
\end{enumerate}
At the limit $[n\to+\infty]$ in \refe{eq:kineticueta}, we obtain \refe{eq:kineticfpre}, so $f$ is a generalized solution to \refe{stoSCL} with initial datum $\mathbf{1}_{u_0>\xi}$.

\subsubsection{Conclusion: proof of Theorem~\ref{th:cvexists}}\label{sec:proofcvexist}
By Theorem~\ref{th:Uadd}, there corresponds a solution $u$ to this $f$: $f=\mathbf{1}_{u>\xi}$. This proves the existence of a solution $u$ to \refe{stoSCL}, unique by Theorem~\ref{th:Uadd}. Besides, owing to the particular structure of $f^\eta$ and $f$, we have
\begin{equation*}
\|u^{\eta_n}\|_{L^2(\T^N\times(0,T))}^2-\|u\|_{L^2(\T^N\times(0,T))}^2=\int_0^T\int_{\T^N}\int_\R 2\xi(f^{\eta_n}-f)d\xi dx dt
\end{equation*}
and (using the bound on $u^\eta$ in $L^3(\T^N)$)
\begin{equation*}
\E\int_0^T\int_{\T^N}\int_{|\xi|>R} |2\xi(f^{\eta_n}-f)(\xi)|d\xi dx dt \leq \frac{C}{1+R}.
\end{equation*}
It follows that $u^{\eta_n}$ converges in norm to $u$ in the Hilbert space $L^2(\T^N\times(0,T)\times\Omega)$. Using the weak convergence, we deduce the strong convergence. Since $u$ is unique, the whole sequence actually converges. This gives the result of the theorem for $p=2$. The case of general $p$ follows from the bound on $u^\eta$ in $L^q$ for arbitrary $q$ and H\"older Inequality. \qed

\begin{appendix}
\section{Proof of Theorem~\ref{th:youngmeasure} and Corollary~\ref{cor:kineticfunctions}}\label{app:proof2}

%

Let $h\in L^1(X)$ be non-negative. By the condition at infinity \refe{nuvanishn}, the sequence of measure $(\nu^n_h)$ defined by
\begin{equation*}
\int_\R \phi(\xi)d\nu^n_h(\xi)=\int_X h(z) \int_\R \phi(\xi)d\nu^n_z(\xi) dz,\quad\phi\in C_b(\R)
\end{equation*} 
is tight. By Prokhorov Theorem, there exists a subsequence still denoted $(\nu^n_h)$ that converges weakly in the sense of measure to a measure $\nu_h$ on $\R$ having the same mass as the measures $(\nu^n_h)$:
\begin{equation}
\nu_h(\R)=\int_X h(z) dz.
\label{totalmass}\end{equation}
Since $L^1(X)$ is separable and $h\mapsto\nu^n_h$ is uniformly continuous in the sense that 
\begin{equation*}
|\nu^n_{h^\flat}(\phi)-\nu^n_{h^\sharp}(\phi)|\leq \|h^\flat-h^\sharp\|_{L^1(X)}\|\phi\|_{C_b(\R)}
\end{equation*}
for all $\phi\in C_b(\R)$, standard diagonal and limiting arguments give $\nu^n_h\rightharpoonup \nu_h$ along a subsequence independent on the choice of $h\in L^1(X)$.
At fixed $\phi\in C_b(\R)$, the estimate $0\leq \nu_h(\phi)\leq \|h\|_{L^1(X)}\|\phi\|_{C_b(\R)}$ and the linearity of $h\mapsto\nu_h(\phi)$ show that
\begin{equation*}
\nu_h(\phi)=\int_X h(z) g(z,\phi)dz,\quad \|g(\cdot,\phi)\|_{L^\infty(X)}\leq\|\phi\|_{C_b(\R)}.
\end{equation*}
Besides, $g(\cdot,\phi)\geq 0$ a.e. since $\nu_h(\phi)\geq 0$ for non-negative $h$, and $\phi\mapsto g(\cdot,\phi)$ is linear. Consequently, for a.e. $z\in X$, we have
\begin{equation*}
g(z,\phi)=\int_\R \phi(\xi) d\nu_z(\phi)
\end{equation*}
where $\nu_z$ is non-negative finite measure on $\R$. By \refe{totalmass}, $\nu_z(\R)=1$. At last, $\nu$ vanishes at infinity since 
\begin{equation*}
\int_X \int_\R |\xi|^p d\nu_z(\xi)d\lambda(z) \leq\limsup_{n\to+\infty}\int_X \int_\R |\xi|^p d\nu^n_z(\xi)d\lambda(z)<+\infty.
\end{equation*}
This concludes the proof of the Theorem. To prove Corollary~\ref{cor:kineticfunctions} we start from the weak convergence $\nu^n_h\rightharpoonup \nu_h$ ($h\in L^1(X)$ being fixed). This implies 
\begin{equation*}
\nu^n_h(\xi,+\infty)\to\nu_h(\xi,+\infty)
\end{equation*}
at all points $\xi$ except the atomic points of $\nu_h$, that are at most countable, hence of zero measure. It follows in particular that
\begin{equation*}
\int_\R \nu^n_h(\xi,+\infty) g(\xi)d\xi\to\int_\R\nu_h(\xi,+\infty)g(\xi)d\xi
\end{equation*}
for all $g\in L^1(\R)$. In other words, we have
\begin{equation}
\int_{X\times\R} f_n(z,\xi)H(z,\xi) d\xi d\lambda(z) \to \int_{X\times\R} f(z,\xi)H(z,\xi) d\xi d\lambda(z)
\label{eq:cvfnH}\end{equation}
for all $H\in L^1(X\times\R)$ of the form $H(z,\xi)=h(z)g(\xi)$. This implies the result since tensor functions are dense in $L^1(X\times\R)$.
\qed

\section{Proof of Lemma~\ref{lem:weakstrongeq}}\label{app:proof3}

By choosing $\theta'(\xi)\gamma(z)$ as a test function, and by use of a standard approximation procedure, we obtain
\begin{equation}
\int_X \theta(u_n(z))\gamma(z) dz\to\int_X\theta(u(z)) \gamma(z) dz
\label{thetagamma}\end{equation}
for all $\theta\in C^1(\R)$, $\gamma\in L^r(X)$ such that $\sup_n\|\theta(u_n)\|_{L^{r'}(X)}<+\infty$, where $r'$ is the conjugate exponent to $r$.
If $p>2$ now, we first obtain the strong convergence in $L^q(X)$, $q=2$, by developing the scalar product
$$
\|u-u_n\|_{L^2(X)}^2=\|u\|_{L^2(X)}^2+\|u_n\|_{L^2(X)}^2-2\<u,u_n\>_{L^2(X)}.
$$
The convergence of the norms follows from \refe{thetagamma} with $\theta(\xi)=\xi$, $\gamma(z)=1$. The weak convergence $\<u,u_n\>_{L^2(X)}\to\|u\|_{L^2(X)}^2$ follows from \refe{thetagamma} with $\theta=1$, $\gamma=u$. Still when $p> 2$, the remaining cases $1\leq q<p$ are obtained by interpolation and by the uniform bound on $\|u_n\|_{L^p(X)}$. If $p\leq2$ now, we notice that, for every $R>0$, the truncate functions 
$$
T_R(u_n):=\min(R,\max(-R,u_n))
$$
satisfy \refe{thetagamma}, and we can apply the reasoning above to show $T_R(u_n)\to T_R(u)$ in $L^r(X)$ strong for every $r<+\infty$. For $1\leq q<p$, the uniform estimate
$$
\|T_R(u_n)-u_n\|_{L^q(X)}\leq \frac{1}{R^{1/s}}\sup_n\|u_n\|_{L^p(X)}^{1+1/s},\quad \frac{1}{q}=\frac{1}{p}+\frac{1}{s},
$$
and similarly for $u$ then gives the result.
\end{appendix}

\providecommand{\bysame}{\leavevmode\hbox to3em{\hrulefill}\thinspace}
\providecommand{\MR}{\relax\ifhmode\unskip\space\fi MR }
\providecommand{\MRhref}[2]{%
  \href{http://www.ams.org/mathscinet-getitem?mr=#1}{#2}
}
\providecommand{\href}[2]{#2}

\end{document}